\documentclass[leqno,oneside,english]{article}
\usepackage[latin1]{inputenc}
\usepackage{amsmath}
\usepackage{amssymb}
\usepackage{color}
\numberwithin{equation}{section}

\makeatletter
\usepackage[T1]{fontenc}
\usepackage[latin1]{inputenc}

\usepackage{graphicx}

\makeatletter \textwidth15cm \textheight20truecm
\newcommand{\eps}{\epsilon}
\newcommand{\R}{{\mathbb R}}

\newcommand{\N}{{\mathbb N}}

\newcommand{\caH}{{\mathcal H}}
\newcommand{\caV}{{\mathcal V}}
\newcommand{\caE}{{\mathcal E}}
\newcommand{\caA}{{\mathcal A}}

\newcommand{\caL}{{\mathcal L}}

\def\be{\begin{equation}}
\def\ee{\end{equation}}
\def\beq{\begin{eqnarray}}
\def\eeq{\end{eqnarray}}
\def\beqs{\begin{eqnarray*}}
\def\eeqs{\end{eqnarray*}}
\def\ea{\end{array}}
\def\ea{\end{array}}
\def\bt{\begin{thm}}
\def\et{\end{thm}}
\def\br{\begin{rk}}
\def\er{\end{rk}}
\def\bl{\begin{lem}}
\def\el{\end{lem}}
\def\bc{\begin{cor}}
\def\ec{\end{cor}}
\def\bex{\begin{exo}}
\def\eex{\end{exo}}
\renewcommand{\div}{\mathop{\rm div}\nolimits}
\newcommand{\rot}{\mathop{\rm curl}\nolimits}

\newtheorem{thm}{Theorem}[section]

\newtheorem{cor}[thm]{Corollary}

\newtheorem{exo}[thm]{Example}

\newtheorem{lem}[thm]{Lemma}

\newtheorem{rk}[thm]{Remark}

\newenvironment{proof}[1][Proof]{\textbf{#1.} }{\ \rule{0.5em}{0.5em}}
\newcommand{\ddiv}{{\rm div}\ }
\newcommand{\Om}{\Omega}

\makeatother

\usepackage{babel}
\makeatother
\author{Serge Nicaise\footnote{Universit\'e Polytechnique Hauts-de-France, C\'ERAMATHS/DEMAV and FR CNRS 2037,
F-59313 - Valenciennes Cedex 9 France,
Serge.Nicaise@uphf.fr}}

\begin{document}
\title{Stability properties of  dissipative evolution equations   with  nonautonomous and nonlinear damping}

  \maketitle

\begin{abstract}
In this paper, we obtain some  stability results of (abstract) dissipative evolution equations   with  a nonautonomous and nonlinear damping using the exponential stability of the retrograde problem
with a linear and autonomous feedback and  a comparison principle. We then  illustrate our abstract statements   for 
different concrete examples, where new results are achieved.
 In a preliminary step, we prove some well-posedness results for some  nonlinear and nonautonomous evolution equations.
\end{abstract}



\section{Introduction}

Stability of evolution equations
of hyperbolic type  with linear or nonlinear autonomous feedbacks
has been  the object of many works.
Let us quote the stability of   the wave equation
\cite{Komornik91,Komornik93b,Komornik93a,Komornikbook,KomornikZuazua,Lasiecka-Triggiani92,zuazua:sicon:90}, 
of the elastodynamic system \cite{AlabauKomornik,Bey:03,guesmia98b,Guesmia99b,guesmia:00,Horn,LaSiam83,Tebou:96},
of the
Petrovsky system \cite{Guesmia98a,Komornikbook,komornik:94}, 
of Mawxell's system
\cite{BH,ELN,KPan,Elleretco2002,Phung} or combination of them 
\cite{daLuz-PerlaMenzala:11,KapiRaupp,nicEES}, see 
also the references cited in the aforementioned works.
On the contrary the case of nonautonomous damping is less considered in the literature, let us quote  
\cite{Daoutli:11,Jiao-Xiao:15,Luo-Xiao:20a,Luo-Xiao:20b,Luo-Xiao:21,martinez:00,mustapha:15,Nakao:97,Pucci-Serrin:96} for the wave equation and  
\cite{BchatniaDaoulatli:13,Bellassoued:08} for the Lam\'e systems.

In the nonautonomous case, even if some similarities appear in the long time behavior of the solution,  the proof is always made  for each particular examples. Hence, our main idea is to  treat the stability of (abstract) evolution equations
of hyperbolic type  with nonautonomous and nonlinear damping by adapting  an approach that was successfully used in the autonomous case in \cite{nicaise:rendiconti:03,Elleretco2002}, namely use
Liu's principle and a comparison principle 
that goes back to \cite{Lasiecka-Tataru:93} and was improved in \cite{Daoutli:11}.
 Liu's  principle 
consists in estimating
the energy of the direct system
by some terms related to the feedbacks
using a retrograde system with final data
equal to the final data of the direct system.
These   terms are then estimated using   the 
  exponential stability of the inverse (retrograde) problem
with a linear feedback (based on Russell's principle)
and a comparison principle. This principle consists in estimating the energy of the systems by the solution of a nonlinear and nonautonomous ODE.
Furthermore, our goal is to present an abstract setting
leading to the stability
of the
abstract (non linear and non autonous) system as soon as the retrograde linear and autonomous system is exponentially stable. Our setting is chosen as large as possible to 
include all examples of the aforementioned papers concerning nonautonomous damping
and allowing  new applications.
The strength of our approach lies in the fact that the  stability results
(with general feedbacks) are only based on the exponential stability  of the retrograde   system 
with a linear  and autonomous feedback,
property that may be checked 
for an explicit problem
by different techniques, like the multiplier method,
microlocal analysis or any method 
entering in a linear framework (like nonharmonic analysis for instance).
We further illustrate our approach by considering
different examples
for which new stability  results
are obtained. Many other examples, like the Petrovsky system or the thermoelastic system, may be treated using the exponential stability  of the retrograde   system 
with a linear  and autonomous feedback, we do not present them for the sake of shortness.

Let us notice that existence results for evolution equations
of hyperbolic type  with nonlinear and nonautonomous feedbacks
are no fully direct, because the domain of the operator may depend on the time variable.
Hence, in a preliminary step, we prove a well-posedness result for a class of  nonlinear and nonautonomous evolution equations,
 extending a result from \cite{Kato67}  and then specializes it to evolution equations 
 of hyperbolic type.


The paper is organized as follows: in section \ref{wellgeneral}
we  give a well posedness result for nonlinear and nonautonomous evolution equations.  In Section \ref{well}, we use this result to obtain some well posedness results for  nonlinear and nonautonomous evolution equations
  of hyperbolic type. Section  \ref{exppol} is devoted to the stability results    for 
a class  of nonautonomous and nonlinear feedbacks
adapting Liu's principle.
Finally in section \ref{sex}  different illustrative examples are treated.


 \section{Well-posedness of nonlinear nonautonomous evolution equations\label{wellgeneral}}
\hspace{5 mm}
\setcounter{equation}{0}

All examples that we will present below can be reduced to a nonlinear evolution equation in a   Hilbert space $X$ of  the form
\begin{equation*}\label{3.1bis}
\left\{
\begin{tabular}{ll}
& $\frac{d U}{dt} (t)
 +A(t)U(t)=0,$  in $X$,\\
& $U(0) = U_0,$
\end{tabular}
\right.
\end{equation*}
where $U$ is the unknown, $U_0\in X$ and $A(t)$ is a (single-valued) nonlinear operators on $X$.
A general theory of such equations with linear operators $A(t)$
has been developed 
using semigroup theory in \cite{KatoCIME,KatoPisa,Pazy} for instance. 
For nonlinear operators $A(t)$
similar results   exist but 
for maximal  quasi-monotone  operators $A(t)$ (for one inner product independent of $t$),
 see \cite{Kato67,Crandal72,Evans77,Lin}
or for maximal  monotone  operators $A(t)$ for a time-dependent
 inner product depending ``smoothly" on $t$, see \cite{NicaisePignotti_03}.
For our systems we need a variant of such results
for
maximal quasi-monotone  operators $A(t)$ for a time-dependent
 inner product depending ``smoothly" on $t$ (see Remarks 4 and 5
in  \cite{Kato67}).
More precisely 
the next result holds.
\bt\label{texinonlinear}
Let $X$  be a  Hilbert space.
For a fixed $T>0$ and any $t\in [0,T]$ we assume that there exists
an inner product $(\cdot,\cdot)_t$ on $X$
depending ``smoothly" on $t$ in the following sense:
there exists $c>0$ such that 
\begin{equation}\label{iptimedependent}
\frac{d}{dt} (u,u)_t\leq 2c (u,u)_t, \forall u\in X, t\in [0,T].
\end{equation}
Furthermore, assume  that 
\\
(i) For all $t\in [0,T]$, $A(t)$ is  single-valued and is a	 maximal quasi-monotone operator
for the inner product $(\cdot,\cdot)_t$, in other words, there exists a non negative real number $\omega$ (independent of $t\in [0,T]$) such that $A(t)+\omega \mathbb{I}$ is a maximal  monotone operator for the inner product $(\cdot,\cdot)_t$,
\\
(ii) the domain  $D(A(t))=D$ of $A(t)$ is independent of $t$, for all $t\in [0,T]$,
\\
(iii) there exists a positive constant $L$
such that
\begin{equation}\label{estderA}
\|A(t)u-A(s)u\|\leq L |t-s|(1+ \|u\|+\|A(s)u\|), \forall u\in D, s,t\in [0,T],
\end{equation}
where for shortness $\|\cdot\|_0$ is denoted by $\|\cdot\|$.
Then for all $a\in D$ the Cauchy problem
\begin{equation}\label{abstractwave}
\left\{
\begin{tabular}{ll}
& $\frac{d u}{dt} (t)
 +A(t)u(t)=0,\ $ for $0\le t\le T$,\\
& $u(0) = a,$
\end{tabular}
\right.
\end{equation}
has a unique solution $u\in C([0,T];X)$ such that $u(t)$
belongs to $D$ for all $t\in [0,T]$,
 its strong derivative 
$\frac{d u}{dt} (t)=-A(t)u(t)$ exists
and is continuous except at a countable number of points $t$.
\et

Note that the condition (\ref{iptimedependent}) and  Gronwall's inequality
imply that
\begin{equation}\label{varnorm}
\|u\|_t\le e^{c|t-s|} \|u\|_s, \forall u\in X, s,t \in [0,T].
\end{equation}
This estimate implies in particular that the norms $\|\cdot\|_t$
are equivalent and gives the  variation of
the norm $\|\cdot\|_t$
with respect to $t$.

\br
{\rm
In the linear case the conditions (\ref{iptimedependent})
and (i) to (iii) imply   that the triplet
$\{A,X,D\}$ forms a CD-system in the sense of \cite{KatoCIME,KatoPisa}.
}
\er

\begin{proof}
The proof is fully similar to the one in  \cite{Kato67}; so, we only give its main steps.
First we recall that $A(t)+\omega \mathbb{I}$ is a  monotone operator for the inner product $(\cdot,\cdot)_t$ if and only if
\be\label{A(t)quasimonotone}
\Re (A(t)u-A(t) v+ \omega(u-v) , u-v)_t\geq 0, \forall u,v\in D,
\ee
or equivalently (see \cite[Lemma 1.1]{Kato67})
\[
\|(1+\alpha \omega)(u-v)+\alpha(A(t)u-A(t)v)_t\geq \|u-v\|_t, \forall u,v\in D,
\alpha>0.
\]
By dividing this estimate by $1+\alpha \omega$ and setting $\lambda=\frac{\alpha}{1+\alpha \omega}$ (that is clearly $<\omega^{-1}$ if $\omega>0$), this is equivalent to
\be\label{A(t)quasimonotoneequiv}
\|u-v+\lambda(A(t)u-A(t)v)_t\geq (1-\lambda \omega) \|u-v\|_t, \forall u,v\in D,
\lambda>0 \hbox{ such that } \lambda \omega<1.
\ee
Hence, we can apply Lemmas 1.1 and 1.2 of \cite{Crandal72}  to $A(t)$  for the norm $t$. In particular for all $n\in \mathbb{N}$
 such that $n>  \omega,$  $\mathbb{I}+n^{-1} A(t)$ is invertible and if we set
\[
J_n(t)=(\mathbb{I}+n^{-1} A(t))^{-1}, \quad A_n(t)= A(t)J_n(t), \forall n\in \mathbb{N} 
\hbox{ such that } n>  \omega,
\]
then the following estimates hold
\beqs
\|J_n(t)x-J_n(t) y\|_t\leq (1-n^{-1}\omega)^{-1} \|x- y\|_t, \forall x,y \in X,
\\
\|A_n(t)x-A_n(t) y\|_t\leq n (1+(1-n^{-1}\omega)^{-1}) \|x- y\|_t, \forall x,y \in X,
\\
\|A_n(t)x\|_t\leq (1-n^{-1}\omega)^{-1} \| A(t)x\|_t, \forall x\in D.
\eeqs
Using \eqref{varnorm}, they are equivalent to
\beq\label{2.4a}
\|J_n(t)x-J_n(t) y\|\leq (1-n^{-1}\omega)^{-1}e^{2cT} \|x- y\|, \forall x,y \in X,
\\
\label{2.4b}
\|A_n(t)x-A_n(t) y\|\leq n (1+(1-n^{-1}\omega)^{-1}) e^{2cT}\|x- y\|, \forall x,y \in X,
\\
\label{2.5}
\|A_n(t)x\|\leq (1-n^{-1}\omega)^{-1} e^{2cT}\| A(t)x\|, \forall x\in D,
\eeq
that respectively correspond to the estimates (2.4) and (2.5) of \cite{Kato67} and are valid for all $n\in \mathbb{N}$
 such that $n>  \omega$. As the factor $(1-n^{-1}\omega)^{-1}e^{2cT}$ is uniformly bounded in $n$ as $n$ goes to infinity, Lemmas 2.4 and 2.5 from \cite{Kato67} remain valid.
Furthermore, by the estimate \eqref{2.5} and our assumption \eqref{estderA}, we have
(see the proof of Lemma 4.1 from \cite{Kato67})
\be
\label{4.2}
\|A_n(t)x-A_n(s)x\|\leq  (1-n^{-1}\omega)^{-1} e^{2cT} L |t-s|(1+ \|u\|+(1+n^{-1})\|A_n(s)u\|), \forall u\in D, s,t\in [0,T], n>\omega,
\ee
that corresponds to the estimate (4.2) of \cite{Kato67}. Since $D$ is dense in $X$, this estimate shows that $A_n(t)$ is Lipschitz continuous in $t$ for all $x\in X$, while \eqref{2.4b} means that the map
$x\to A(t) x$ is Lipschitz continuous for a fixed $t\in [0,T]$, uniformly in $x$ and $t$. Thus the approximated problem
 \begin{equation}\label{abstractwaven}
\left\{
\begin{tabular}{ll}
& $\frac{du_n}{dt} (t)
 +A_n(t)u_n(t)=0,\ $ for $0\le t\le T$,\\
& $u_n(0) = a,$
\end{tabular}
\right.
\end{equation}
has a unique solution $u_n\in C^1([0,T];X)$ for all $a\in X$.
We now show that the statements of Lemma 4.2 of \cite{Kato67} hold if $a\in D$, namely there exists a positive constant $K$ (that  depends on $c$,   $\omega$,  $T$,  and $\|a\|+\|A(0) a\|$ but not on $n$) such that
\beq\label{Le4.2/1}
\|u_n(t)\|\leq K, \forall t\in [0,T], n>\omega,
\\
\label{Le4.2/2}
 \|u'_n(t)\|=\|A_n(t)u_n(t)\|\leq K, \forall t\in [0,T], n>\omega,
 \eeq
 where for shortness we write $\frac{du_n}{dt}=u'_n$.
 Indeed for $t\in [0,T)$, let us fix $h$   in $[0, T-t]$ and  set $x_n(t):= u_n(t+h)-u_n(t)$. As $x_n$ is differentiable in $t$ and usng \eqref{iptimedependent}, we have
 \[
2\|x_n(t)\|_t   \frac{d}{dt} \|x_n(t)\|_t=\frac{d}{dt} \|x_n(t)\|_t^2\leq 2c \|x_n(t)\|_t^2+ 2\Re (x_n'(t), x_n(t))_t.
 \]
Using  \eqref{abstractwaven}, we get
\beqs
\|x_n(t)\|_t   \frac{d}{dt} \|x_n(t)\|_t &\leq &c \|x_n(t)\|_t^2-\Re(A_n(t+h) u_n(t+h)-A_n(t) u_n(t), x_n(t))_t
 \\
 &\leq &c \|x_n(t)\|_t^2-\Re(A_n(t+h) u_n(t+h)-A_n(t+h) u_n(t), x_n(t))_t
 \\
 &-&
 \Re(A_n(t+h) u_n(t)-A_n(t) u_n(t), x_n(t))_t.
 \eeqs
 Using \eqref{A(t)quasimonotone} and \eqref{4.2}, we obtain
 \beqs
\|x_n(t)\|_t   \frac{d}{dt} \|x_n(t)\|_t 
 &\leq &(c+\omega) \|x_n(t)\|_t^2 
 \\
 &+&
  (1-n^{-1}\omega)^{-1} e^{2cT} L h(1+ \|u_n(t)\|+(1+n^{-1})\|u'_n(t)\|) \|x_n(t)\|_t,
 \eeqs
 Simplifying by $\|x_n(t)\|_t$ (see \cite[p. 515]{Kato67}), we find
  \beqs
  \frac{d}{dt} \|x_n(t)\|_t 
 \leq (c+\omega) \|x_n(t)\|_t 
 +
  (1-n^{-1}\omega)^{-1} e^{2cT} L h(1+ \|u_n(t)\|+(1+n^{-1})\|u'_n(t)\|).
 \eeqs
 This estimate directly implies that
 \[
  \frac{d}{dt} \left(e^{-(c+\omega) t} \|x_n(t)\|_t\right) \leq 
    L_1 h(1+ \|u_n(t)\|+\|u'_n(t)\|),
    \]
    for a positive constant $L_1$ that depends on  $c$,   $\omega$ and $T$ but is independent of $n$. Integrating this estimate in $(0,t)$,  we find
    \[
   e^{-(c+\omega) t} \|x_n(t)\|_t-\|x_n(0)\|_0 \leq 
    L_1 h\int_0^t (1+ \|u_n(s)\|+ \|u'_n(s)\|)\,ds.
    \]
 By  \eqref{varnorm}, we find
  \[
   \|x_n(t)\|\leq  L_2(\|x_n(0)\|+
    h\int_0^t (1+ \|u_n(s)\|+ \|u'_n(s)\|)\,ds),
    \]
    for a positive constant $L_2$ that depends on  $c$,   $\omega$ and $T$ but is independent of $n$.
  Dividing by $h$ and letting $h$ goes to zero, we obtain
  \[
   \|u'_n(t)\| \leq  L_2(\|u'_n(0)\| +
    \int_0^t (1+ \|u_n(s)\|+ \|u'_n(s)\|)\,ds).
    \]
  As $u'_n(0)=A(0) a$ and 
  \[
  u_n(t)=a+\int_0^t u'_n(s)\,ds,
  \]
  we find as in \cite[p. 516]{Kato67} that
  \[
  \|u_n(t)\|+ \|u'_n(t)\| \leq  L_3 (1+
    \int_0^t (\|u_n(s)\|+ \|u'_n(s)\|)\,ds),
    \]
for   a positive constant $L_3$ that depends on  $c$,   $\omega$,  $T$ and 
$\|a\|+\|A(0) a\|$ but is independent of $n$. By Gronwall's Lemma, we deduce that 
  \eqref{Le4.2/1} and \eqref{Le4.2/2} hold.
  
 We now show that   the statements of Lemma 4.3 of \cite{Kato67} hold, namely for $a\in D$, 
 the strong limit $u(t)=\lim_{n\to \infty} u_n(t)$ exists uniformly for $t\in [0,T]$
 and $u$ is Lipschitz continuous. Indeed for all $m,  n\in \mathbb{N}$ such that $m, n>\omega$, we set
 $x_{mn}(t)=u_m(t)-u_n(t)$ and as before we have
 \[
 \frac{d}{dt} \|x_{mn}(t)\|_t^2\leq 2c \|x_{mn}(t)\|_t^2+ 2\Re (x_{mn}'(t), x_{mn}(t))_t.
 \]
 Using \eqref{abstractwaven} and  \eqref{A(t)quasimonotone}, we find
 \[
 \frac{d}{dt} \|x_{mn}(t)\|_t^2\leq 2c  \|x_{mn}(t)\|_t^2+ 2\omega \|y_{mn}(t)\|_t^2
 +2\Re (A_m(t) x_m(t)-A_n(t) x_n(t), y_{mn}(t)-  x_{mn}(t))_t,
 \]
 where $y_{mn}(t)=J_m(t)u_m(t)-J_n(t)u_n(t)$. By the triangle inequality, we then have
 \[
 \frac{d}{dt} \|x_{mn}(t)\|_t^2\leq 2(c+2\omega)  \|x_{mn}(t)\|_t^2+ 4\omega \|y_{mn}(t)-x_{mn}(t)\|_t^2+
 2\Re (A_m x_m(t)-A_n(t) x_n(t), y_{mn}(t)-  x_{mn}(t))_t,
 \]
 Using the estimate \eqref{Le4.2/2} and \eqref{varnorm}, we arrive at
 \[
 \frac{d}{dt} \|x_{mn}(t)\|_t^2\leq 2(c+2\omega)  \|x_{mn}(t)\|_t^2+ K_1 \|y_{mn}(t)-x_{mn}(t)\|^2+
  K_1 \|y_{mn}(t)-  x_{mn}(t)\|,
 \]
for  a positive constant $K_1$ that  depends on $c$,   $\omega$,  $T$,  and $\|a\|+\|A(0) a\|$ but not on $m, n$. Obviously, this is equivalent to
 \[
 \frac{d}{dt}\left(e^{-2 (c+2\omega)t} \|x_{mn}(t)\|_t^2\right)\leq   K_1(\|y_{mn}(t)-x_{mn}(t)\|^2+
   \|y_{mn}(t)-  x_{mn}(t)\|),
 \]
 and integrating it between $(0,t)$, we find (as $x_{mn}(0)=0$)
 \[
 e^{-2 (c+2\omega)t} \|x_{mn}(t)\|_t^2\leq K_1\int_0^t(\|y_{mn}(s)-x_{mn}(s)\|^2+
 \|y_{mn}(s)-  x_{mn}(s)\|)\,ds.
  \]
  This finally leads to
  \[
 \|x_{mn}(t)\|^2\leq e^{2(3c+2\omega)T} K_1 \int_0^t( \|y_{mn}(s)-x_{mn}(s)\|^2+
  \|y_{mn}(s)-  x_{mn}(s)\|\,ds.
  \]
  As 
  \[
  y_{mn}(s)-  x_{mn}(s)= J_m(s) u_m(s)-u_m(s)+J_n(s) u_n(s)-u_n(s)=n^{-1}A_n(s) u_n(s)-m^{-1}A_m(s) u_m(s),\] by \eqref{Le4.2/2}, we obtain
  \[
 \| y_{mn}(s)-  x_{mn}(s)\|\leq K (m^{-1}+n^{-1}).
 \]
 Inserting this estimate in the previous one, we arrive at
 \[
 \|x_{mn}(t)\|^2\leq K_2 (m^{-1}+n^{-1}), \forall t\in [0,T],
 \]
 for  a positive constant $K_2$ that  depends on $c$,   $\omega$,  $T$,  and $\|a\|+\|A(0) a\|$ but not on $m, n$. Thus the strong limit $u(t)=\lim_{n\to \infty} u_n(t)$ exists uniformly in $t\in [0,T]$. The Lipschitz continuity of $u$ follows from the uniform Lipschitz property of the $u_n$, that is consequence of \eqref{Le4.2/2}.
  
 The remainder  of the proof is the same as in \cite{Kato67} since it is based on the properties proved before.
\end{proof}

\br{\rm
Obviously, the previous Theorem remains valid if $X$ is a real Hilbert space.
}
\er

\section{Abstract hyperbolic setting\label{well}}
In this section we describe a general abstract setting 
 of hyperbolic
type inspired from \cite{nicaise:rendiconti:03}
that will be used later on.
It is motivated by
 the examples (and other ones) given in section \ref{sex}
which all enter in this setting.

\subsection{General assumptions\label{ssgeneralhypo}}

Let us fix two real  Hilbert spaces $\caH$,
$\caV$ with respective inner products
$(.,.)_\caH$, $(.,.)_\caV$
and such that 
$\caV$ is densely and continuously embedded into $\caH$.
Identifying $\caH$ with its dual $\caH'$ we have the standard
diagram
\[
\caV\hookrightarrow \caH=\caH'\hookrightarrow \caV'.
\]
The duality pairing between $\caV'$ and $\caV$ will be denoted by $\langle \cdot,\cdot \rangle$, so that 
\[
\langle u,v\rangle= (u,v)_\caH, \forall u,v\in \caH.
\]
We suppose that $\caV$ is continuously embedded into   a  control space $U$, that is supposed to be 
in the form
\be\label{defU}
U=\prod_{j=1}^J U_j,
\ee
where for all $j=1,\cdots, J\in \N^\star:=\N\setminus\{0\}$,
$U_j$ is a closed subspace of 
$L^2(X_j,\mu_j)^{N_j}$, with $N_j\in \N^\star$, 
when $X_j$ is a metric space, and $(X_j,{\cal A}_j,\mu_j)$ is a measure space such that $\mu_j(X_j)<\infty$.

  For all $j=1,\cdots, J$, we suppose given a mapping $\alpha_j\in C([0,\infty)\times X_j; (0,\infty))$ and  locally Lipschitz with respect to the time variable, in the sense that 
  for all $T$, there exist a positive constant $\kappa(T)$ (that may depend on $T$) such that
  \be\label{lipintime}
 |\alpha_j(t,x)- \alpha_j(t,x)|\leq  \kappa(T) |t-s|, \forall t \in
[0,T], x\in X_j,
\ee
 and  a continuous mapping
$g_j:\R^{N_j}\to \R^{N_j}$ such that
\beq
\label{gmono}
&& (g_j(x)-g_j(y))\cdot (x-y)\geq 0, \forall x,y\in  \R^{N_j}\hbox{ (monotonicity)},
\\
\label{goo}
&& g_j(0)=0,
\\
\label{gbounded}
&& |g_j(x)|\leq M(1+|x|), \forall x\in  \R^{N_j}, 
\eeq
for some positive constant $M$.  

We further define  the (nonlinear) time-dependent operator $B(t)$ from $\caV$ into $\caV'$   by
\be\label{defBnonlineaire}
\langle B(t)u,v\rangle=\sum_{j=1}^J\int_{X_j} \alpha_j(t,x_j) g_j(( I_Uu)_j(x_j))\cdot ( I_Uv)_j(x_j)\,d\mu_j(x_j), \forall u,v\in \caV,
\ee
where  $I_U$ denotes the embedding from $\caV$ to $U$ 
and therefore, $( I_Uu)_j$ is the $j^{th}$ component of $I_Uu$.

We finally suppose given a bounded linear operator $A_1$
from $\caV$ into $\caV'$ and consider the evolution equation
\be\label{a6}
\left\{
\begin{tabular}{lllll}
& $\frac{dx}{d t}(t)+ A_1 x(t)+B(t)x(t)=0 $ in $\caH, t\geq 0$,\vspace{2mm}
\\
& $x(0) =  x_0$.
\end{tabular}
\right.
\ee
 
 This system clearly involves the  (nonlinear) and time-dependent operator
$\caA_B(t)$  defined by
\beq\label{a2}
&&D(\caA_B(t))=\{v\in \caV|(A_1+ B(t))v\in\caH\},
\\
\label{a3}
&&\caA_B(t)=(A_1+B(t))v, \forall v\in D(\caA_B(t)).
\eeq
In its full generality, the domain of $\caA_B(t)$ depends on the time variable. Consequently we cannot apply Theorem \ref{texinonlinear}. Nevertheless there are two cases treated below for which this Theorem applies. In both cases, if $x_0\in D(\caA_B(0))$, we will show that a unique solution $x$ exists with the following properties:
\be\label{regsolutiona6}
\left\{
\begin{tabular}{ll}
$x\in   C([0,\infty), \caH)$ is such that $x(t)\in D(\caA_B(t)),$ for all $t\in [0,\infty)$
 and \\
 $x'(t)=-\caA_B(t) x(t)$ exists in $\caH$
and is continuous except at a countable number of points $t$.
\end{tabular}
\right.
\ee

Before going on let us show that under the additional assumption
 that
\be 
\label{h4} \langle A_1 u,u\rangle=0, \forall u\in \caV,
\ee
system \eqref{a6} is dissipative.

\bl\label{l2.1}
Under the above assumptions,
for all $t\geq 0$, the operator $A_B(t)$ is monotone for the natural inner product of $\caH$. namely
\be\label{AB(t)mononotone}
(A_B(t) u-A_B(t) v,u-v)_\caH=\langle B(t) u-B(t) v, u-v\rangle\geq 0, \forall u\in D(\caA_B(t)).
\ee
Consequently if $x$ is a solution of \eqref{a6} with the regularity \eqref{regsolutiona6},
  its associated energy
\be\label{a7}
\caE(t)=\frac12||x(t)||_\caH^2
\ee
is non-increasing; moreover,
we have
\beq\label{2.2weak}
&&\caE(S)-\caE(T)=\int_S^T\langle B(t)u(t),u(t)\rangle\,dt, \forall 0\leq S<T<\infty,
\\
\label{a9}
&&\frac{d }{dt}\caE(t)=-\langle B(t)u(t),u(t)\rangle \leq 0, \hbox{ for a. a. } t\geq 0.
\eeq
\el
\begin{proof}
Let us first show that $A_B(t)$ is monotone. Indeed for any $u,v \in D(\caA_B(t))$, by the
definition of $A_B(t)$  and  the property  \eqref{h4}, we have
\beqs
&&(A_B(t) u-A_B(t)v, u-v)_\caH=\langle A_1(u-v), u-v\rangle
\\
&&\hspace{1cm}+\langle B(t) u-B(t) v, u-v\rangle
=\langle B(t) u-B(t) v, u-v\rangle.
\eeqs
Finally by the 
definition of $B(t)$ and then \eqref{gmono} and recalling that $\alpha_j(t,x)>0$, we have
\beqs
&&\langle B(t) u-B(t) v, u-v\rangle
\\
&&\hspace{1cm}=
\sum_{j=1}^J\int_{X_j} \alpha_j(t,x_j)\left(g_j(( I_Uu)_j(x_j))-g_j(( I_Uv)_j(x_j))\right)\cdot \left((I_Uv)_j(x_j)-(I_Uv)_j(x_j)\right)\,d\mu_j(x_j)
\\
&&\hspace{1cm}\geq  0,
\eeqs
which proves \eqref{AB(t)mononotone}.
 
For the second assertion
it suffices to show (\ref{a9})
since  \eqref{2.2weak} follows by integration between $S$ and $T$.
By the regularity assumptions on $x$, we have
\[
\frac{d }{dt}\caE(t)=(x'(t),x(t))_\caH=-(\caA_B(t) u(t),u(t))_\caH, \hbox{ for a. a.} t\geq 0,
\]
by (\ref{a6}).
By our assumption \eqref{goo}, we have $A_B(t) 0=0$
and consequently by \eqref{AB(t)mononotone}, we get  (\ref{a9}).
\end{proof}

\subsection{The  ``bounded''  case}

We here assume that $\caH$ is continuously embedded into $U$. As we shall see below this assumption implies that $B(t)$ becomes a (nonlinear) operator from $\caH$ into itself and therefore, the domain of $\caA_B(t)$ does not depend on $t$ anymore. 

\bt\label{texistencebdB}
In addition to the previous assumptions, assume that $\caH$ is continuously embedded into $U$,
and that there exists a positive real number $\lambda$ such that the range $R(\lambda \mathbb{I}+\caA_B(t))$ is equal to $\caH$. Then \be\label{Dt=D}
D(\caA_B(t))=D=\{u\in \caV| A_1 u\in\caH\}, \forall t\geq 0,
\ee
and
for any $x_0\in D,$ problem \eqref{a6}
has a unique solution  $x$ satisfying \eqref{regsolutiona6}.
\et
\begin{proof}
We first show that
\eqref{Dt=D} holds.
Indeed as  $\caH$ is continuously embedded into $U$, the mapping $I_U$ extends to a linear and continuous operator from $\caH$ into $U$; therefore, there exists a positive constant $C$ such that
\be\label{IUbd}
\|I_U u\|_\caH\leq C \|u\|_\caH, \forall u\in \caH.
\ee
By our assumption \eqref{gbounded} and the definition of $B(t)$, we then have
\[
|\langle B(t)u,v\rangle|\leq M\sum_{j=1}^J\int_{X_j}  \alpha_j(t,x_j) (1+
|(I_Uu)_j(x_j)|) | I_Uv)_j(x_j|\,d\mu_j(x_j), \forall u,v\in \caV.
\]
By the continuity property of $\alpha_j(t,\cdot)$, Cauchy-Schwarz's inequality
and the estimate \eqref{IUbd}, we obtain
\[
|\langle B(t)u,v\rangle|\leq C(t) (1+\|u\|_\caH)\|v\|_\caH, \forall u,v\in \caV,
\]
where $C(t) $ is a positive constant that depends on $M, C$,  and $t$.
As $\caV$ is dense in $\caH$, for a fixed $u\in \caV$, the linear mapping
\[
\caV\to \R: v\to \langle B(t)u,v\rangle,
\]
can be extended to a linear and continuous form to the whole $\caH$. 
By the Riesz's representation theorem, there exists $h(t)\in \caH$ such that
\[
\langle B(t)u,v\rangle=(h(t), v)_\caH, \forall v\in \caH.
\]
In other words, for $u\in \caV$, $B(t)u$ can be identified with $h(t)$ and therefore,
$(A_1+ B(t))u\in\caH$ if and only if $A_1 u\in\caH$, which proves \eqref{Dt=D}.

By Lemma \ref{l2.1} and our additional assumption $R(\lambda \mathbb{I}+\caA_B(t))=\caH$, for some $\lambda>0$, we deduce that the assumption (i) of 
Theorem \ref{texinonlinear} holds.

Let us end up with the third assumption.  Fix $T>0$ and let   $u \in D$,  and $t,s\in [0,T]$, then  we clearly have  
\[
\caA_B(t) u-\caA_B(s)u=B(t) u-B(s)u.
\]
Therefore, for any $v\in \caH$,  by the definition of $B(t)$ and our previous considerations, we may write
\[
(\caA_B(t) u-\caA_B(s)u, v)_\caH=
\sum_{j=1}^J\int_{X_j} (\alpha_j(t,x_j)-\alpha_j(s,x_j)) g_j(( I_Uu)_j(x_j))\cdot ( I_Uv)_j(x_j)\,d\mu_j(x_j).
\]
By our assumptions \eqref{lipintime} to \eqref{gbounded}, we obtain
\beqs
|(\caA_B(t) u-\caA_B(s)u, v)_\caH|&\leq& \kappa(T) |t-s|
\sum_{j=1}^J\int_{X_j} g_j(( I_Uu)_j(x_j))\cdot ( I_Uv)_j(x_j)\,d\mu_j(x_j)
\\
&\leq& \kappa(T) |t-s|
\sum_{j=1}^J\int_{X_j} (1+
|(I_Uu)_j(x_j)|)  \cdot ( I_Uv)_j(x_j)\,d\mu_j(x_j).
\eeqs
 Cauchy-Schwarz's inequality
and the estimate \eqref{IUbd} allow to conclude that
\[
|(\caA_B(t) u-\caA_B(s)u, v)_\caH|\leq \sqrt{2} C \kappa(T) |t-s| 
(\sum_{j=1}^J \mu_j(X_j)+C \|u\|_\caH) \|v\|_\caH).
\]
Since this estimate is valid for all $v\in \caH$, this means that
\[
\|\caA_B(t) u-\caA_B(s)u\|_\caH\leq \sqrt{2} C \kappa(T) |t-s| 
(\sum_{j=1}^J \mu_j(X_j)+C \|u\|_\caH),
\]
and proves that the assumption (iii) of 
Theorem \ref{texinonlinear} holds.  

In conclusion by Theorem \ref{texinonlinear} for $x_0\in D$ and  any $T>0$, there exists a unique solution $u_T\in C([0,T];\caH)$ of problem 
\be\label{a6local}
\left\{
\begin{tabular}{lllll}
& $\frac{dx_T}{d t}(t)+ A_1 x_T(t)+B(t)x_T(t)=0 $ in $\caH, t\in [0,T]$,\vspace{2mm}
\\
& $x_T(0) =  x_0,$
\end{tabular}
\right.
\ee
such that $x_T(t)$
belongs to $D$ for all $t\in [0,T]$,
 its strong derivative 
$\frac{d x_T}{dt} (t)=-A(t)x_T(t)$ exists
and is continuous except at a countable number of points $t$.
By  uniqueness, for $T'>T$, the restriction of $x_{T'}$ to $[0,T]$ coincides with $x_T$. Therefore, a unique global solution   
  $x\in C([0,\infty);\caH)$ of \eqref{a6} exists with the properties  \eqref{regsolutiona6}.
\end{proof}

%
 
 \subsection{The  ``unbounded''  case\label{abunbdcase}}
 
 Here we assume that the mappings
 $\alpha_j$ do  not depend on the $x_j$ variable and coincide, namely
 there exists a mapping $\alpha\in C^1([0,\infty; (0,\infty))$ 
such that $\alpha'$ is   locally Lipschitz (in the sense that 
for all $T>0$, there exists a positive constant $\nu(T)$ such that
$
|\alpha'(t)-\alpha'(s)|\leq \nu(T) |t-s|,$ for all $s,t\in [0,T]
$)
 such that
 \be\label{alphaequal}
\alpha_j(t, x_j)=\alpha(t), \forall x_j\in X_j, t\geq 0.
\ee
Due to \eqref{defBnonlineaire}, this means that 
$B(t)=\alpha(t) B_1$, where
\be\label{defBnonlineairepartcase}
\langle B_1u,v\rangle=\sum_{j=1}^J\int_{X_j} g_j(( I_Uu)_j(x_j))\cdot ( I_Uv)_j(x_j)\,d\mu_j(x_j), \forall u,v\in \caV.
\ee

\bt\label{texistencebubdB}
In addition to the   assumptions made in subsection \ref{ssgeneralhypo}, assume that \eqref{alphaequal} holds, that $A_1+B_1$ is maximal quasi monotone with a dense domain in $\caH$, and that there exist two mappings 
$D\in C^1([0,\infty);\caL(\mathcal{\caH}))$ and $\tilde D\in C([0,\infty);\caL(\mathcal{\caH}))$
such that $D'$ and $\tilde D$ are locally Lipschitz and for all $t\geq 0$,
$D(t)$ and $\tilde D(t)$ are invertible, $D(t) \tilde D(t)$ is symmetric positive definite and for all $T>0$,
there exists a positve constant $c_T$ such that
\be\label{equivnorm}
(\tilde D(t)^{-1} D(t) ^{-1}x,x)_\caH\geq c_T\|x\|_\caH^2, \forall x\in \caH, \forall t\in [0,T],
\ee
and finally
\be\label{commutateur}
(A_1+\alpha(t) B_1) D(t)^{-1}= \tilde D(t) (A_1+ B_1), \forall t\geq 0.
\ee
Then for all $x_0\in D(\caA_B(0))$,  problem \eqref{a6}
has a unique solution  $x$ satisfying \eqref{regsolutiona6}.
\et
\begin{proof}
Assuming that the solution $x$ of  problem \eqref{a6} exists and is smooth enough, we perform the change of unknown
\[
\tilde x(t)=D(t) x(t).
\]
Hence, as $\tilde x'(t)=D'(t) x(t)+D(t) x'(t)$ and by \eqref{a6} , we get
\[
\tilde x'(t)=D'(t) D(t)^{-1} \tilde x(t)-D(t) (A_1+\alpha(t) B_1) D(t)^{-1} \tilde x(t).
\]
With our assumption \eqref{commutateur}, we arrive at
\be\label{17/9:1}
\tilde x'(t)=D'(t) D(t)^{-1} \tilde x(t)-D(t) \tilde D(t) (A_1+ B_1)\tilde x(t).
\ee
This corresponds to \eqref{abstractwave} with the operator 
\[
A(t) =D(t) \tilde D(t) (A_1+ B_1)-D'(t) D(t)^{-1}, 
\]
whose domain is clearly
\[
D(A(t))=D(A_1+B_1), 
\]
and is independent of $t$, due to our assumptions on $D(t)$ and  $\tilde D(t)$.

In order to apply  Theorem \ref{texinonlinear} we introduce the time dependent inner product
\[
(x, \tilde x)_t=( \tilde D(t)^{-1}D(t)^{-1} x, \tilde x)_\caH, \forall x, \tilde x\in \caH.
\]
Our assumptions on $D$ and $\tilde D$ guarantee that it is indeed an inner product on $\caH$ whose associated  norm  is equivalent to the standard one, namely for a fixed $T$, we have
\be\label{equivnormtinnerprod}
\sqrt{c_T} \|x\|_\caH\leq \|x\|_t\leq  C_T \|x\|_\caH, \forall x\in \caH,
\ee
for some positive constant   $C_T$, and  that the property \eqref{iptimedependent} also holds. From its definition, we see that
$A(t)$ is quasi monotone for this inner product. Indeed from it definition,
for any $x,y\in D(A_1+B_1)$, and $t\in [0,T]$, we have
\[
(A(t) x-A(t) y, x-y)_t=((A_1+ B_1) x-(A_1+ B_1) y,x-y)_\caH-(\tilde D(t)^{-1}D(t)^{-1} D'(t) D(t)^{-1} (x-y),x-y)_\caH.
\]
Hence, as $A_1+B_1$ is quasi monotone in $\caH$ (i.e. $A_1+B_1+\omega_1 \mathbb{I}$ is monotone for some $\omega_1\geq 0$), and due to our assumptions on $D$ and $\tilde D$, we then have
\[
(A(t) x-A(t) y, x-y)_t\geq -\omega_T \|x-y\|^2_\caH,
\]
for some $\omega_T>0$ (depending on $T$). Due to the equivalence \eqref{equivnormtinnerprod}, we arrive at 
\[
(A(t) x-A(t) y, x-y)_t\geq -\omega_T C_T^2 \|x-y\|_t^2,
\]
which yields the  quasi monotonicity of $A(t)$. Let us now show the maximality property.
Indeed for $\lambda>0$ large enough, we want to show that
$E(t):=\lambda \mathbb{I}+D(t) \tilde D(t) (A_1+ B_1)-D'(t) D(t)^{-1}$ is surjective.
But as $D(t) \tilde D(t)$ is an isomorphism, this is equivalent to the surjectivity of
$E(t):=\lambda (D(t) \tilde D(t))^{-1}+ A_1+ B_1-(D(t) \tilde D(t))^{-1}D'(t) D(t)^{-1}$.
Now we take advantage of Theorem 1 of
\cite{Browder:68} by considering the previous operator as a perturbation of 
$T_1:=A_1+B_1+\omega_1 \mathbb{I}$ (that satisfies the assumption of this Theorem). Due to the linearity of $T_2(t):=\lambda (D(t) \tilde D(t))^{-1}-(D(t) \tilde D(t))^{-1}D'(t) D(t)^{-1}-\omega_1 \mathbb{I}$, it is clearly hemicontinous
and due to the assumption \eqref{equivnorm}, for $\lambda>0$ large enough, 
$T_2(t)$ will be monotone, bounded and coercive. Using the above Theorem, we deduce that $E(t)=T_1+T_2(t)$
is surjective.
In summary assumption (i) of Theorem \ref{texinonlinear} holds and it remains to check
the assumption (iii) of this Theorem. For that purpose, let us fix  $x\in D(A_1+B_1)$
and $s,t\in [0,T]$, then by definition we have
\[
A(t) x-A(s) x=(D(t) \tilde D(t)-D(s) \tilde D(s)) (A_1+ B_1)x
+(D'(s) D(s)^{-1}-D'(t) D(t)^{-1})x
\]
By the local Lipschitz property of $D$, $\tilde D$ and of  the derivative of $D$, we get
\[
\|A(t) x-A(s) x\|_\caH\leq K(T) |t-s| (\|(A_1+ B_1)x\|_\caH+\|x\|_\caH).
\]
We now transform  
\[
(A_1+ B_1)x=(D(s) \tilde D(s))^{-1}
 (D(s) \tilde D(s) (A_1+ B_1)x-D'(s) D(s)^{-1})+(D(s) \tilde D(s))^{-1}D'(s) D(s)^{-1}x,
 \]
 use the triangle inequality and use the continuity of $D$, $\tilde D$ and $D'$ to find 
 \[
\|A(t) x-A(s) x\|_\caH\leq K_1(T) |t-s| (\|A(s) x\|_\caH+\|x\|_\caH).
\]
for  a positive constant $K_1(T)$,
 which implies that \eqref{estderA} is valid.
 
 In conclusion by Theorem \ref{texinonlinear}, there exists a unique solution 
 $\tilde x$ of  \eqref{17/9:1} with initial condition $\tilde x(0)=D(0)x_0$ (that belongs to $D(A_1+B_1)$ by the assumption on $x_0$)
 satisfying 
 $\tilde x\in   C([0,\infty), \caH)$,  $x(t)\in D(A_1+B_1),$ for all $t\in [0,\infty)$
 and 
 $x'(t)=-A(t) x(t)$ exists in $\caH$
and is continuous except at a countable number of points $t$.
Setting $x(t)=D(t)^{-1} \tilde x(t)$, we readily check that it is the unique solution of problem \eqref{a6} and that it satisfies \eqref{regsolutiona6}.
\end{proof}

\section{Stability results in the nonlinear and nonautonomous case \label{exppol}}

Here we use Liu's principle  \cite{Liu}  
and a comparison principle with a nonlinear and nonautonomous ODE  from  \cite{Daoutli:11}  (see also \cite{Lasiecka-Tataru:93}) to deduce  decay  rates of the energy 
using appropriate nonlinear and nonautonomous feedbacks.

We first  recall the  comparison principle  obtained in \cite{Daoutli:11}  (compare with    \cite[Theorem 2 and Corollary 2]{Lasiecka-Tataru:93})).
\bt\label{tii}
Let $\beta$ be a continuous mapping from $[0,\infty)$ to $(0,\infty)$ and $p$ a strictly increasing convex  mapping  
from $[0,+\infty)$ to $[0,+\infty)$ such that $p(0)=0$.
Let $\caE:[0,+\infty)\to [0,+\infty)$ be a non-increasing mapping satisfying
\be\label{A6}
  \beta((n+1)T) p(\caE(nT))+\caE((n+1)T) \leq \caE(nT), \forall n\in \mathbb{N},
\ee
for some $T>0$.  
%
Then 
\be\label{A7}
\caE(t)\leq \psi^{-1}\left(\int_T^t \beta(s)\, ds\right), \forall t\geq T,
\ee
where $\psi$  is defined by
\be\label{A13,5}
\psi(x)=\int_x^{\caE(0)}\frac{1}{p(s)}\,ds, \forall x> 0.
\ee
\et
\begin{proof}
Let us shortly recall the proof from \cite{Daoutli:11}. 
Since \eqref{A7} trivially holds if $\caE(0)=0$ (because in such a case $\caE(t)=0$, for all $t\geq 0$), we can assume that $\caE(0)>0$. First by \cite[Lemma 4.2]{Daoutli:11}, 
 the next comparison principle holds
\be\label{comparison}
\caE(t)\leq S(t-T), \forall t\geq T,
\ee
where $S$ is the unique solution of the nonlinear and nonautonomous ODE 
\be\label{ODE}
S'(t)+\beta(t+T) p(S(t))=0, \forall t\geq 0, \quad S(0)=\caE(0).
\ee
Such a solution exists and remains positive for all $t>0$ due the Cauchy-Lipschitz Theorem (because the assumptions on $p$ garantee that it is locally Lipschitz in $[0,\infty)$). 

With the help of \cite[Lemma 4.2/2.]{Daoutli:11} (the properties on $p$ guarantee that the assumption (24) from \cite{Daoutli:11}  holds with $m=p(1)^{-1}$), we deduce that 
\[
S(t)\leq \psi^{-1}\left(\int_0^t \beta(s+T)\, ds\right), \forall t\geq 0,
\]
with $\psi$ defined by \eqref{A13,5} (and is meaningful   because 
$\lim_{x\to 0+}\psi(x)=+\infty$ 
reminding that $p(x)\leq p(1) x$, for all $x\in [0,1]$). This estimate combined with \eqref{comparison} yields the result.
\end{proof}

Let us now  recall Russell's principle that yields an  exact controllability result for the evolution equation 
associated with the operator $-A_1$ with controls in $L^2(]0,T[;U)$
provided $A_1-I_U$ generates  a semigroup of contractions and 
$-A_1-I_U$ generates 
an exponentially stable semigroup  of contractions in $\caH$, see
\cite[Theorem 4.1]{nicaise:rendiconti:03}.
\bt\label{tec}
Assume that $A_1-I_U$ 
generates  a semigroup of contractions  in $\caH$ and that $-A_1-I_U$
generates a  semigroup of contractions $S(t)$ in $\caH$ that is   exponentially stable in the sense that there exists two positive constants $C$ and $\omega$ such that
\be\label{expdecay}\|S(t) x_0\|_\caH\leq C e^{-\omega t} \| x_0\|_\caH, \forall x_0\in \caH.
\ee
Then there exists $T>0$ large enough, such that for any $p_0\in\caH$, there exists a control $K\in L^2((0,T);U)$  
 such that the solution  $p\in C([0,T]; \caH)$ of
\be\label{ce1inv}
\left\{
\begin{tabular}{lllll}
& $\frac{\partial p}{\partial t} +A_1 p=K $ in $\caV', t\in [0, T]$,\vspace{2mm}
\\
& $p(T) =  p_0,$
\end{tabular}
\right.
\ee
satisfies
\be\label{ce2inv}
p(0)=0.
\ee
Furthermore, there exists a positive constant $D>1$ depending only on $T$, and the constants $C$ and $\omega$ such that
\be\label{spo4}
\int_0^T \|K(t)\|_U^2\,dt+\int_0^T \|I_Up(t)\|_U^2\,dt
 \leq
2D  \|p_0\|_\caH^2.
\ee
\et

We now give the consequence of this result to our  system (\ref{a6}) in three different  cases of functions $\alpha_j$: non-increasing with respect to $t$, non-decreasing with respect to $t$ and oscillating  with respect to $t$. But first we give an energy estimate valid in all cases.

\bl\label{lestET}
Under the assumptions of Theorem \ref{tec},  any solution $x$ of 
(\ref{a6}) with the regularity \eqref{regsolutiona6}, 
satisfies
\be\label{spo9bis} 
\caE(T)\leq D
\left(\sum_{j=1}^J\int_0^T 
\int_{X_j}
\{|( I_Ux(t))_j(x_j)|^2+\alpha_j(t,x_j)^2
|g_j(( I_Ux(t))_j(x_j))|^2\}\,d\mu_j(x_j) 
\,dt\right).
\ee
\el
\begin{proof}
Let $x$  be the unique
solution of (\ref{a6}) satisfying \eqref{regsolutiona6} 
and consider $p$ the solution of
problem (\ref{ce1inv}) and (\ref{ce2inv}) with $p_0=x(T)\in \caH$  with
$T>0$ from  Theorem
\ref{tec}.  
Owing to \cite[Remark 4.2]{nicaise:rendiconti:03}, consider   a sequence $p_\epsilon\in W^{1,\infty}([0,\infty), \caH)\cap L^\infty([0,\infty), \caV)$ of strong solution of 
(\ref{ce1inv}) with final data $p_{0\epsilon}$  tending to $p$ in
 $C([0,T],\caH)$ as $\epsilon $ goes to zero and satisfying 
 \beq
&&\label{RKeps} K_\epsilon\to K \hbox{ in } L^2(]0,T[;U) \hbox{ as } \epsilon\to 0,
\\
&&\label{Rpeps} I_U p_\epsilon\to I_U  p \hbox{ in } L^2(]0,T[;U) \hbox{ as } \epsilon\to 0.
\eeq
 By (\ref{a6}) and (\ref{ce1inv})  we may write
 $$
\langle\partial_t x+A_1x+B(t)x,p_\epsilon\rangle_{\caV',\caV}+
 \langle\partial_t
p_\epsilon+A_1p_\epsilon-K_\epsilon,x\rangle_{\caV',\caV}=0, \hbox{ for a.a. } t\in [0,T].
 $$
As the assumption (\ref{h4}) yields
$$
\langle A_1x,p_\epsilon\rangle_{\caV',\caV}+
\langle A_1p_\epsilon,x\rangle_{\caV',\caV}=0,
$$
the above identity reduces to
\beqs
(\partial_t x,p_\epsilon)_\caH+(\partial_t p_\epsilon,x)_\caH
+
\langle B(t)x,p_\epsilon\rangle_{\caV',\caV}-
\langle K_\epsilon,x\rangle_{\caV',\caV}
=0,  \hbox{ for a.a. } t\in [0,T].
\eeqs
Integrating this identity for $t\in(0,T)$, we get
$$
(x(T),p_\epsilon(T))_\caH-
(x(0),p_\epsilon(0))_\caH+\int_0^T
(\langle B(t)x,p_\epsilon\rangle_{\caV',\caV}-
\langle K_\epsilon,x \rangle_{\caV',\caV})\,dt=0.
$$
By the definitions of $K_\epsilon$ 
and $B(t)$ we arrive at 
\beqs
(x(T),p_\epsilon(T))_\caH- 
(x(0),p_\epsilon(0))_\caH&=&
\int_0^T 
\Big((K_\epsilon, I_Ux)_{U} 
\\
&-&\sum_{j=1}^J\int_{X_j}\alpha_j(t,x_j) g_j(( I_Ux)_j(x_j))\cdot ( I_Up_\epsilon)_j(x_j)\,d\mu_j(x_j) 
\Big)\,dt.
\eeqs
Passing to the limit in $\epsilon$ and 
using the  initial and final conditions  on $p$, we  
obtain
\[
2\caE(T)= \int_0^T
\left((K, I_Ux)_{U}
-\sum_{j=1}^J\int_{X_j}\alpha_j(t,x_j)  g_j(( I_Ux)_j(x_j))\cdot ( I_Up)_j(x_j)\,d\mu_j(x_j)
\right)\,dt.
\]
Cauchy-Schwarz's inequality   leads finally to
\beq\label{spo9}
2\caE(T)&\leq &\|K\|_{L^2(0,T;U)}\| I_Ux\|_{L^2(0,T;U)}
\\
\nonumber
&+& \| I_Up\|_{L^2(0,T;U)}
\left(\sum_{j=1}^J\int_0^T  
\int_{X_j} \alpha_j(t,x_j)^2 |g_j(( I_Ux)_j(x_j))|^2\,d\mu_j(x_j)
\,dt\right)^{1/2}.
\eeq
Using the estimate \eqref{spo4}   (recalling that $p_0=x(T)$), we have
\[
\int_0^T \|K(t)\|_U^2\,dt+\int_0^T \|I_Up(t)\|_U^2\,dt
 \leq 4D
  \caE(T).
\]
Using this   estimate in the previous one,  we arrive at
\eqref{spo9bis}.
\end{proof}

\bc\label{cestEn+1T}
Under the assumptions of Theorem \ref{tec},  any solution $x$ of 
(\ref{a6}) with the regularity \eqref{regsolutiona6}, 
satisfies
\be\label{spo9bisn+1} 
\caE((n+1)T)\leq D
\Big(\sum_{j=1}^J\int_{nT}^{(n+1)T} 
\int_{X_j}
\{|( I_Ux(t))_j(x_j)|^2+\alpha_j(t,x_j)^2
|g_j(( I_Ux(t))_j(x_j))|^2\}\,d\mu_j(x_j) 
\,dt\Big),
\ee
for all $n\in \N$.
\ec
\begin{proof}
We apply the previous Lemma to $x_n$ (instead of $x$)
defined by
\[
x_n(t)=x(t+nT), \forall t\geq 0,
\]
that is still solution of (\ref{a6}) with the regularity \eqref{regsolutiona6}, where the (nonlinear) and time-dependent operator $B$ is replaced by
$B_n(t)=B(t+nT)$. The estimate \eqref{spo9bis} applied to $x_n$
 yields 
 \[
\caE((n+1)T)\leq D
\Big(\sum_{j=1}^J\int_{0}^{T}\!\!
\int_{X_j}\!
\{|( I_Ux(t+nT))_j(x_j)|^2+\alpha_j(t+nT,x_j)^2
|g_j(( I_Ux(t+nT))_j(x_j))|^2\}\,d\mu_j(x_j) 
\,dt\Big),
\]
 that is nothing else than \eqref{spo9bisn+1} by a simple change of variable.
\end{proof}

\subsection{The non-increasing case}
\bt\label{tstabobsgen}
In addition to the previous assumptions on $g_j$ and $\alpha_j, j=1,\cdots, J,$ suppose that
  $g_j$  satisfies   
\beq
\label{gstrp}
&& g_j(x)\cdot x\geq m |x|^2, \forall x\in  \R^{N_j}: |x|\geq 1,
\\
\label{gnear0}
&&|x|^2+|g_j(x)|^{2}\leq G( g_j(x)\cdot x), \forall  x\in  \R^{N_j}: |x|\leq 1,
\eeq
for some positive constant $m$ and a concave strictly increasing function $G:[0,\infty)\to [0,\infty)$
such that $G(0)=0$. Furthermore, we assume that
for all $j=1,\cdots, J$ and all $x_j\in X_j,$
 the mapping
\be\label{ajdecroissante}
\alpha_j(\cdot, x_j): [0,\infty)\to (0,\infty): t\to \alpha_j(t, x_j)
\hbox{ is   non-increasing, }
\ee
\be\label{upperbdalphaj}
\tilde \alpha:=\max_{1\leq j\leq J} \sup_{x_j\in X_j} \alpha_j(0,x_j) <\infty,
\ee
and 
\be\label{alphaequiv} 
\alpha(t)=\min_{1\leq j\leq J} \inf_{x_j\in X_j} \alpha_j(t,x_j)>0, \forall t\in [0,\infty).
\ee

Under the assumptions of Theorem \ref{tec}, there exists   $c>0$  
 (depending on $T$, $C$, $\omega$ (from Theorem \ref{tec}), $\max_j\mu_j(X_j)$, $\tilde \alpha$,    $M$, and  $m$) such that
\be\label{A19}
\caE(t)\leq \psi^{-1}\left(T\mu \int_T^t \alpha(s)\, ds\right), \forall t\geq T,
\ee
for all solution $x$ of 
(\ref{a6}) satisfying \eqref{regsolutiona6}, 
 where $\mu=\min_j\mu_j(X_j)$, $\psi$ is given by
(\ref{A13,5}) with $p=h^{-1}$, and
$h$ defined by
\be\label{defh}
h(x)=c (x+G(x)), \forall x\geq 0,
\ee
\et
\begin{proof}
Let $x$  be the unique
solution of (\ref{a6}) satisfying \eqref{regsolutiona6}
and let $n$ be an arbitrary nonnegative integer.
Using \eqref{spo9bisn+1} and the definition of $\tilde\alpha$, we get
\be\label{spo9bisalphadecroissant}
\caE((n+1)T)\leq C_1
\Big(\sum_{j=1}^J\int_{nT}^{(n+1)T}
\int_{X_j}
\{|( I_Ux)_j(x_j)|^2+\alpha_j(t,x_j)
|g_j(( I_Ux)_j(x_j))|^2\}\,d\mu_j(x_j) 
\,dt\Big),
\ee
with $C_1=D \max\{1,\tilde \alpha\}$.

We now estimate   the right-hand side of (\ref{spo9bisalphadecroissant}) as follows:
For all $j=1,\cdots, J$,
introduce 
\beq
\label{s+}
&& \Sigma_{j,n}^+=\{ (x,t)\in X_j\times (nT,(n+1)T)\ :\ |(I_U x)_j(x,t)|>1\},
\\
&& \Sigma_{j,n}^-=\{ (x,t)\in X_j\times (nT,(n+1)T) \ :\ |(I_U x)_j(x,t)|\leq 1\},
\label{s-}
\eeq
and   split up
\[\int_{nT}^{(n+1)T}
\int_{X_j}\{|( I_Ux)_j(x_j)|^2+\alpha_j(t,x_j)
|g_j(( I_Ux)_j(x_j))|^2\}\,d\mu_j(x_j)\,dt=I_{j,n}^++I_{j,n}^-,
\]
where
\beqs
I_{j,n}^+&:=&
\int_{\Sigma_{j,n}^+}\{|( I_Ux)_j(x_j)|^2+\alpha_j(t,x_j)
|g_j(( I_Ux)_j(x_j))|^2\}\,d\mu_j(x_j)\,dt,
\\
I_{j,n}^-&:=&
\int_{\Sigma_{j,n}^-}\{|( I_Ux)_j(x_j)|^2+\alpha_j(t,x_j)
|g_j(( I_Ux)_j(x_j))|^2\}\,d\mu_j(x_j)\,dt.
\eeqs

For the estimation of $I_{j,n}^+$, we first  notice that
the assumption  (\ref{gbounded})  leads to
\[
I_{j,n}^+\leq 
\int_{\Sigma_{j,n}^+}(1+2M\alpha_j(t,x_j))|( I_Ux)_j(x_j)|^2\,d\mu_j(x_j)\,dt,
\]
and  by (\ref{gstrp})  and \eqref{upperbdalphaj} we get
\[
I_{j,n}^+\leq   m^{-1} (1+ 2M \tilde \alpha) \int_{\Sigma_{j,n}^+}  ( I_Ux)_j(x_j)\cdot g_j(( I_Ux)_j(x_j)) \,d\mu_j(x_j)dt. 
\]
As $\alpha_j(\cdot, x_j)$ is non-increasing, and using
\eqref{upperbdalphaj}-\eqref{alphaequiv}, we have
\be\label{13/9:1}
1\leq \frac{ \alpha_j(t, x_j)}{\alpha_j((n+1)T, x_j)}\leq \frac{ \alpha_j(t, x_j)}{\alpha((n+1)T)},
\forall x_j\in X_j, t\in [nT,(n+1)T],
\ee
which  allows to obtain 
\be\label{estI^+}
I_{j,n}^+\leq   m^{-1} (1+ 2M \tilde \alpha) \alpha((n+1)T)^{-1} \int_{\Sigma_{j,n}^+}   \alpha_j(t, x_j) ( I_Ux)_j(x_j)\cdot g_j(( I_Ux)_j(x_j)) \,d\mu_j(x_j)dt.
\ee
Since  (\ref{gmono}) and (\ref{goo}) yield
\be\label{gpos}
g_j(x)\cdot x\geq 0, \forall x\in \R^{N_j},
\ee
and since $ \alpha_j(t, x_j) > 0$ for all $t$ and $x_j\in X_j$,  
we have
\beq\label{serge:11/10:1}
 &&\int_{\Sigma_{j,n}^+}  \alpha_j(t, x_j) ( I_Ux)_j(x_j)\cdot g_j(( I_Ux)_j(x_j)) \,d\mu_j(x_j)dt
\\
&&\hspace{-1cm}
\leq
\int_{nT}^{(n+1)T}
\int_{X_j}  \alpha_j(t, x_j) ( I_Ux)_j(x_j)\cdot g_j(( I_Ux)_j(x_j)) \,d\mu_j(x_j)dt
\leq  (\caE(nT)-\caE((n+1)T)),
\nonumber
\eeq
this last estimate following from \eqref{2.2weak}.
Using this estimate in \eqref{estI^+},
we arrive at
\be\label{spo10}
I_{j,n}^+
\leq c_1 \alpha((n+1)T)^{-1}   (\caE(nT)-\caE((n+1)T)),
\ee
for some positive constant $c_1$ depending only on $\tilde \alpha$,  $M$ and $m$.

Similarly by the assumption (\ref{gnear0})
and the monotonicity of $G$ and $\alpha$ we have 
\beqs
I_{j,n}^-&\leq& \max\{1,\tilde \alpha\}
\int_{\Sigma_{j,n}^-}G(( I_Ux)_j(x_j)\cdot g_j(( I_Ux)_j(x_j))) \,d\mu_j(x_j)dt
\\
&\leq&  \max\{1,\tilde \alpha\} \int_{nT}^{(n+1)T}
\int_{X_j} G(( I_Ux)_j(x_j)\cdot g_j(( I_Ux)_j(x_j))) \,d\mu_j(x_j)dt.
\eeqs
Jensen's inequality then yields
\[
I_{j,n}^-\leq   \max\{1,\tilde \alpha\} T \mu_j(X_j) 
G\left(\frac{1}{T\mu_j(X_j)}  \int_{nT}^{(n+1)T}
\int_{X_j}  ( I_Ux)_j(x_j)\cdot g_j(( I_Ux)_j(x_j)) \,d\mu_j(x_j)dt
\right).
\]
As  $G$ is strictly increasing and again using \eqref{13/9:1}, we obtain 
 \beqs
I_{j,n}^-\leq   
K
G\left(\frac{1}{T\mu_j(X_j)\alpha((n+1)T)} \!\! \int_{nT}^{(n+1)T}
\int_{X_j} \alpha_j(t,x_j)( I_Ux)_j(x_j)\cdot g_j(( I_Ux)_j(x_j)) \,d\mu_j(x_j)dt
\right),
\eeqs
where  $K= \max\{1,\tilde \alpha\} T \max_j \mu_j(X_j)$.
By  (\ref{2.2weak}), we arrive at
\be\label{spo11}
I_{j,n}^-\leq \ K
G\left(\frac{\caE(nT)-\caE((n+1)T)}{T \mu_j(X_j)\alpha((n+1)T)}
\right).
\ee

The estimates   (\ref{spo10}) and  (\ref{spo11})  into the estimate
(\ref{spo9bisalphadecroissant}) and the monotonicity of $G$ give
 \[
\caE((n+1)T)
\leq c_2    \left\{\frac{\caE(nT)-\caE((n+1)T)}{T\mu \alpha((n+1)T)}
+G\left(\frac{\caE(nT)-\caE((n+1)T)}{T\mu \alpha((n+1)T)}\right)\right\},
\]
for some positive constant $c_2$ (depending on $T$, $C$, $\omega$ (from  Theorem \ref{tec}),  $\max_j\mu_j(X_j)$, $\tilde \alpha$, $M$ and $m$), recalling that $\mu=\min_j\mu_j(X_j)$.
As \eqref{upperbdalphaj}-\eqref{alphaequiv} imply that $\alpha((n+1)T)\leq \tilde \alpha$, this finally leads to
\beqs
\caE(nT)&=&\caE(nT)-\caE((n+1)T)+\caE((n+1)T)
\\
&\leq& \max\{\mu T  \tilde \alpha,c_{2}\}
\left\{\frac{\caE(nT)-\caE((n+1)T)}{T\mu \alpha((n+1)T)}
+G\left(\frac{\caE(nT)-\caE((n+1)T)}{T\mu \alpha((n+1)T)}\right)\right\}.
\eeqs
With $c=\max\{\mu T \tilde \alpha,c_{2}\}$, and the definition \eqref{defh} of $h$, we have found that
\[
\caE(nT)\leq h\left(\frac{\caE(nT)-\caE((n+1)T)}{T\mu \alpha((n+1)T)}\right),
\]
which can be equivalently written as
%
%
%
%
\be\label{A3}
T\mu  \alpha((n+1)T) h^{-1}(\caE(nT))+\caE((n+1)T) \leq \caE(nT).
\ee
Since this estimate is valid for all  $n \in \mathbb{N}$, we conclude by Theorem \ref{tii}
with the choice $\beta(t)=T\mu  \alpha(t)$.
\end{proof}


Note that the conditions \eqref{gbounded} and \eqref{gstrp} means that $g_j$ is linearly bounded at infinity; therefore, the decay rate in (\ref{A19}) is guided by the behaviour of 
$g_j$ near zero and of the behavior of $\int_0^t\alpha(s)\,dx$ as $t$ goes to $\infty$.
Since we are mainly interested in the influence of the time dependency on the decay rate, we restrict ourselves to   examples of functions $g_j$ that are linear, sublinear or superlinear near 0 (compare with subsection 3.2.1 and Example 1 of \cite{Daoutli:11}).

\bex\label{exemple}
{\rm
Suppose that  $g_j$ satisfies  (\ref{gmono}) to (\ref{gbounded})
and (\ref{gstrp}) as well as
\be\label{spo15} 
x\cdot g_j(x)\geq c_0|x|^{p+1}, \,
|g_j(x)|\leq C_0|x|^{\gamma}, \forall |x|\leq 1,
\ee
for some positive constants $c_0, C_0$, $\gamma\in (0,1]$ and $p\geq \gamma$. Then $g_j$ satisfies 
(\ref{gnear0}) with
$G(x)=x^{\frac{2}{q+1}}$ and $q=\frac{p+1}{\gamma}-1$ (which is $\geq 1$).

If $p=\gamma=1$ (then $q=1$), that corresponds to a linear behavior of $g_j$ near 0,
we have $G(x)=x$ and, hence, $h(x)=2c x$. Therefore,  under the other assumptions of  Theorem \ref{tstabobsgen}
we get the
 decay 
 \[
 \caE(t)\leq K  \caE(0) e^{-L\int_0^t\alpha(s)\,ds}, \forall t\geq 0,
 \]
 for some positive constants $K$ and $L$, since $\psi^{-1}(t)=\caE(0)e^{-\frac{t}{2c}}$.
 
On the contrary if $p+1>2\gamma$ (corresponding to the sublinear case if $p=2$
and to the superlinear case if $\gamma=1$ and $p>1$), then  we get  the decay
$K(\caE(0))\left(\int_0^t\alpha(s)\,ds\right)^{-\frac{2\gamma}{p+1-2\gamma}}$ (since $\psi^{-1}(t)$ is equivalent to $t^{\frac{2}{1-q}}$ for $t$ large), with 
$K(\caE(0))= K 
(1+\caE(0)^{-\frac{p+1-2\gamma}{2\gamma}})^{-\frac{2\gamma}{p+1-2\gamma}}$, with a positive constant $K$.

Note that in both cases, the energy tends to zero as soon as
\[
\int_0^t\alpha(s)\,ds\to \infty, \hbox{ as } t \to \infty.
\]
In particular, if $\alpha(t)=\frac{1}{(1+t)^\sigma}$, with $\sigma>0$, in both cases, we get
  \[
 \caE(t)\leq K(\caE(0)) t^{-r},
 \]
 for some $r>0$ (with $K(\caE(0))= K \caE(0)$ in the linear case) that, in the linear case, translates an underdamped  phenomenon.

A function $g$ satisfying all these assumptions is given by
$$
g(x)=
\left\{
\begin{tabular}{ll}
$|x|^{\gamma-1} x$ & if $ |x|\leq 1,$\\
$x$ & if $ |x|\geq 1,$
\end{tabular}
\right.
$$
for some $\gamma\in (0,1]$. In that case (\ref{spo15}) holds for $p=\gamma$.
}
\eex

\subsection{The non-decreasing case}
\bt\label{tstabobsgencroissant}
In addition to the assumptions on $g_j$ and $\alpha_j, j=1,\cdots, J,$ from subsection 
\ref{ssgeneralhypo}, suppose that
  $g_j$  satisfies   \eqref{gstrp} and \eqref{gnear0}  
for some positive constant $m$ and a concave strictly increasing function $G:[0,\infty)\to [0,\infty)$
such that $G(0)=0$ and satisfying the additional assumption
\be\label{propG}
\exists \delta \geq 2, C_G>0: \beta^2 G(x)\leq C_G G(\beta^\delta x), \forall x, \beta\in (0,\infty).
\ee
Furthermore, we assume that
for all $j=1,\cdots, J$ and all $x_j\in X_j,$
 the mapping
\be\label{ajcroissante}
\alpha_j(\cdot, x_j): [0,\infty)\to (0,\infty): t\to \alpha_j(t, x_j)
\hbox{ is   non-decreasing, }
\ee
and that for all $t\in [0,\infty)$
\be\label{alphasup} 
\alpha(t)=\max_{1\leq j\leq J} \sup_{x_j\in X_j} \alpha_j(t,x_j)<\infty,
\ee
and
\be\label{alphazero} \alpha(0)>0,
\ee 
so that the mapping
\[
\alpha: [0,\infty)\to (0,\infty) : t\to \alpha(t)
\]
is non-decreasing. We finally suppose that there exists $c_0\in (0,1]$ such that
\be\label{equivalence} 
c_0\alpha(t)\leq  \alpha_j(t,x_j) \leq \alpha(t), \forall t\in [0,\infty), x_j\in X_j, j=1,\cdots, J.
\ee

Under the assumptions of Theorem \ref{tec}, there exists   $c>0$  
 (depending on $T$, $C$, $\omega$ (from  Theorem \ref{tec}), $\max_j\mu_j(X_j)$, $\alpha(0)$,  $c_0$,  $M$, and  $m$) such that
\be\label{A19croissante}
\caE(t)\leq \psi^{-1}\left(T\mu c_0\int_{T}^{t}  \alpha(s-T) \alpha(s)^{-\delta}\, ds\right), \forall t\geq T,
\ee
for all solution $x$ of 
(\ref{a6}) satisfying \eqref{regsolutiona6}, 
 where $\mu=\min_j\mu_j(X_j)$, $\psi$ is given by
(\ref{A13,5}) with $p=h^{-1}$, and
$h$ is defined by \eqref{defh}.
\et
\begin{proof}
Let $x$  be the unique
solution of (\ref{a6}) satisfying \eqref{regsolutiona6}
and let $n$ be arbitrary nonnegative integer.
We   estimate   the right-hand side of \eqref{spo9bisn+1} as follows: Using the  sets 
$\Sigma_{j,n}^+$ and $\Sigma_{j,n}^-$ defined by \eqref{s+} and \eqref{s-} respectively, we     split up
\[\int_{nT}^{(n+1)T}
\int_{X_j}\{|( I_Ux)_j(x_j)|^2+\alpha_j(t,x_j)^2
|g_j(( I_Ux)_j(x_j))|^2\}\,d\mu_j(x_j)\,dt=I_{j,n}^++I_{j,n}^-,
\]
where
\beqs
I_{j,n}^+&:=&
\int_{\Sigma_{j,n}^+}\{|( I_Ux)_j(x_j)|^2+\alpha_j(t,x_j)^2
|g_j(( I_Ux)_j(x_j))|^2\}\,d\mu_j(x_j)\,dt,
\\
I_{j,n}^-&:=&
\int_{\Sigma_{j,n}^-}\{|( I_Ux)_j(x_j)|^2+\alpha_j(t,x_j)^2
|g_j(( I_Ux)_j(x_j))|^2\}\,d\mu_j(x_j)\,dt.
\eeqs

For the estimation of $I_{j,n}^+$, we first  notice that
the assumptions  \eqref{alphazero} and \eqref{equivalence} 
lead to
\beq\label{serge11/10:2}
\alpha_j(t,x_j)&\leq& \alpha((n+1)T), \forall t\in [nT, (n+1) T],
\\
c_0 \alpha(0)&\leq& \alpha_j(t,x_j), \forall t\geq 0.
\label{serge11/10:3}
\eeq
 Therefore, using the assumption (\ref{gbounded}) on $g_j$, we have
\[
I_{j,n}^+\leq (\frac{1}{c_0^2 \alpha(0)^2}+2M) \alpha((n+1)T)
\int_{\Sigma_{j,n}^+}   \alpha_j(t,x_j)|( I_Ux)_j(x_j)|^2\,d\mu_j(x_j)\,dt,
\]
and  by (\ref{gstrp}) we get
\[
I_{j,n}^+\leq   m^{-1}  (\frac{1}{c_0 \alpha(0)}+2M) \alpha((n+1)T) \int_{\Sigma_{j,n}^+}  \alpha_j(t,x_j)| ( I_Ux)_j(x_j)\cdot g_j(( I_Ux)_j(x_j)) \,d\mu_j(x_j)dt. 
\]
Since the estimate \eqref{serge:11/10:1} remains valid, we obtain
\be\label{spo10croissante}
I_{j,n}^+
\leq c_1 \alpha((n+1)T)   (\caE(nT)-\caE((n+1)T)),
\ee
for some positive constant $c_1$ depending only on $c_0$, $ \alpha(0)$,  $M$ and $m$.

Let us go on with the estimation of $I_{j,n}^-$.
First using \eqref{serge11/10:2}-\eqref{serge11/10:3}, we may write 
\[
I_{j,n}^-\leq C_1 \alpha((n+1)T)^2
\int_{\Sigma_{j,n}^-} \{|( I_Ux)_j(x_j)|^2+ 
|g_j(( I_Ux)_j(x_j))|^2\}\,d\mu_j(x_j)\,dt,
\]
where $C_1=\max\{\frac{1}{c_0^2\alpha(0)^2}, 1\}$.
Hence, by the assumption (\ref{gnear0})
and the monotonicity of $G$ and the positivity of $\alpha_j$, as before we have 
\[
I_{j,n}^- 
\leq C_1\alpha((n+1)T)^2  \int_{nT}^{(n+1)T}
\int_{X_j} \alpha_j(t,x_j) G(( I_Ux)_j(x_j)\cdot g_j(( I_Ux)_j(x_j))) \,d\mu_j(x_j)dt.
\]
Jensen's inequality then yields
\[
I_{j,n}^-\leq   C_1\alpha((n+1)T)^2  T \mu_j(X_j) 
G\left(\frac{1}{T \mu_j(X_j)}  \int_{nT}^{(n+1)T}
\int_{X_j} ( I_Ux)_j(x_j)\cdot g_j(( I_Ux)_j(x_j)) \,d\mu_j(x_j)dt
\right).
\]
Now   \eqref{alphazero} and \eqref{equivalence} yield
\[
c_0\alpha(nT)\leq \alpha_j(t,x_j), \forall t\in [nT, (n+1)T],
\] and since $G$ is strictly increasing, we then obtain 
 \beqs
I_{j,n}^-\leq   
C_2 \alpha((n+1)T)^2
G\left(\frac{1}{T\mu c_0\alpha(nT)} \!\! \int_{nT}^{(n+1)T}
\int_{X_j} \alpha_j(t,x_j)( I_Ux)_j(x_j)\cdot g_j(( I_Ux)_j(x_j)) \,d\mu_j(x_j)dt
\right),
\eeqs
where  $C_2= C_1T\max_j\mu_j(X_j)$.
By  (\ref{2.2weak}), we arrive at
\be\label{spo11croissante}
I_{j,n}^-\leq \ C_2  \alpha((n+1)T)^2
G\left(\frac{\caE(nT)-\caE((n+1)T)}{T\mu c_0\alpha(nT)}
\right).
\ee
At this stage, we take advantage of the property 
\eqref{propG} 
to conclude that
\be\label{spo11croissante2}
I_{j,n}^-\leq \ C_2
G\left(\frac{\alpha((n+1)T)^\delta(\caE(nT)-\caE((n+1)T))}{T\mu c_0\alpha(nT)}
\right).
\ee

The estimates   (\ref{spo10croissante}) (as $\alpha((n+1)T)\leq  \frac{\alpha((n+1)T)^\delta}{ \alpha(0)^{\delta-2}\alpha(nT)}$ because $\alpha$ is non-decreasing and $\delta\geq 2$) and  (\ref{spo11croissante2})  into the estimate
\eqref{spo9bisn+1}  give
 \[
\caE((n+1)T)
\leq c_2    \left\{\frac{\alpha((n+1)T)^\delta(\caE(nT)-\caE((n+1)T))}{T\mu c_0\alpha(nT)}
+G\left(\frac{\alpha((n+1)T)^\delta(\caE(nT)-\caE((n+1)T))}{T\mu c_0\alpha(nT)}\right)\right\},
\]
for some positive constant $c_2$ (depending on $T$,  $\max_j\mu_j(X_j)$, $c_0$, $\alpha(0),$ $\delta$, $C$, $\omega$, $M$ and $m$).
As  the non-decreasing property of $\alpha$ implies that
 $\frac{\alpha((n+1)T)^\delta}{\alpha(nT)}\geq  \alpha(0)^{\delta-1}$, this finally leads to
\beqs
\caE(nT)&=&\caE(nT)-\caE((n+1)T)+\caE((n+1)T)
\\
&\leq& c
\left\{\frac{\alpha((n+1)T)^\delta(\caE(nT)-\caE((n+1)T))}{T\mu c_0\alpha(0)}
+G\left(\frac{\alpha((n+1)T)^\delta(\caE(nT)-\caE((n+1)T))}{T\mu c_0\alpha(0)}\right)\right\},
\eeqs
where $c= \max\{\frac{T\mu c_0}{\alpha(0)^{\delta-2}},c_{2}\}$. By the definition \eqref{defh} of $h$, we have found that
\[
\caE(nT)\leq h\left(\frac{\alpha((n+1)T)^\delta(\caE(nT)-\caE((n+1)T))}{T\mu c_0\alpha(nT)}\right),
\]
which can be equivalently written as
%
%
%
%
\be\label{A3croissante}
T\mu c_0\alpha(nT) \alpha((n+1)T)^{-\delta} h^{-1}(\caE(nT))+\caE((n+1)T) \leq \caE(nT).
\ee
Since this estimate is valid for all  $n \in \mathbb{N}$, we conclude by Theorem \ref{tii}
with the choice $\beta(t)=T\mu   c_0 \alpha(t-T) \alpha(t)^{-\delta}$.
\end{proof}

\bex\label{exemplecroissante}
{\rm
If $g_j$ satisfies the assumptions from Example \ref{exemple}, $G$ is given by
$G(x)=x^{\frac{2}{q+1}}$ and $q=\frac{p+1}{\gamma}-1\geq 1$; hence,
it satisfies the assumption \eqref{propG} with $C_G=1$
and $\delta=q+1=\frac{p+1}{\gamma}$.

If $p=\gamma=1$ (then $q=1$), that corresponds to a linear behavior of $g_j$ near 0,
we have $G(x)=x$ and, hence, $h(x)=2c x$. Under the other assumptions of  Theorem \ref{tstabobsgencroissant}
we then get the
 decay 
 \[
 \caE(t)\leq K  \caE(0) e^{-L\int_T^t\alpha(s-T) \alpha(s)^{-2}\,ds}, \forall t\geq T,
 \]
 for some positive constants $K$ and $L$, since $\psi^{-1}(t)=\caE(0)e^{-\frac{t}{2c}}$.   
 
On the contrary if $p+1>2\gamma$ (corresponding to the sublinear case if $p=2$
and to the superlinear case if $\gamma=1$ and $p>1$), then  we get  the decay
$K(\caE(0))\left(\int_T^t\alpha(s-T) \alpha(s)^{-\frac{p+1}{\gamma}}\,ds\right)^{-\frac{2\gamma}{p+1-2\gamma}}$ (since $\psi^{-1}(t)$ is equivalent to $t^{\frac{2}{1-q}}$ for $t$ large).  

Note that in both cases, the energy tends to zero as soon as
\[
 \int_T^t\alpha(s-T) \alpha(s)^{-\delta}\,ds\to \infty, \hbox{ as } t \to \infty.
\]

In particular, if $\alpha(t)=(1+t)^\sigma$, with $0<\sigma\leq \frac{1}{\delta-1}=\frac{\gamma}{p+1-\gamma}$,  in both cases, we get
  \[
 \caE(t)\leq K(\caE(0))t^{-r},
 \]
 for some $r>0$, that, in the linear case, translates an overdamping phenomenon.
}
\eex
 
 \subsection{The oscillating case}
  
 \bt\label{tstabobsgenosc}
In addition to the assumptions on $g_j$ and $\alpha_j, j=1,\cdots, J,$ from subsection 
\ref{ssgeneralhypo}, suppose that
  $g_j$  satisfies   \eqref{gstrp} and \eqref{gnear0}  
for some positive constant $m$ and a concave strictly increasing function $G:[0,\infty)\to [0,\infty)$
such that $G(0)=0$ Furthermore, we assume that there exists two positive constants $\alpha_0$ and $\tilde \alpha$ such that
\be\label{alphaoscj}
  \alpha_0\leq \alpha_j(t,x_j) \leq \tilde \alpha,  \forall t\in [0,\infty), x_j\in X_j, j=1,\cdots, J.
\ee

Under the assumptions of Theorem \ref{tec}, there exists   $c>0$  
 (depending on $T$, $C$, $\omega$ (from Theorem \ref{tec}), $\max_j\mu_j(X_j)$, $\tilde \alpha$, $\alpha_0$,  $M$, and  $m$) such that
\be\label{A19osc}
\caE(t)\leq \psi^{-1}\left(T\mu \alpha_0 (t-T) \right), \forall t\geq T,
\ee
for all solution $x$ of 
(\ref{a6}) satisfying \eqref{regsolutiona6}, 
 where $\mu=\min_j\mu_j(X_j)$, $\psi$ is given by
(\ref{A13,5}) with $p=h^{-1}$, and
$h$ defined by
\eqref{defh}.
\et
\begin{proof}
Let $x$  be the unique
solution of (\ref{a6}) satisfying \eqref{regsolutiona6}
and let $n$ be arbitrary nonnegative integer.
Using \eqref{spo9bisn+1} and the assumption \eqref{alphaoscj}, we get
\be\label{spo9bisalphaosc}
\caE((n+1)T)\leq C_1
\Big(\sum_{j=1}^J\int_{nT}^{(n+1)T}
\int_{X_j}
\{|( I_Ux)_j(x_j)|^2+
|g_j(( I_Ux)_j(x_j))|^2\}\,d\mu_j(x_j) 
\,dt\Big),
\ee
with $C_1=D \max\{1,\tilde \alpha^2\}$.

We now 
estimate   the right-hand side of (\ref{spo9bisalphaosc}) as follows:  Using the  sets 
$\Sigma_{j,n}^+$ and $\Sigma_{j,n}^-$ from \eqref{s+} and  \eqref{s-}, we     split up
\[\int_{nT}^{(n+1)T}
\int_{X_j}\{|( I_Ux)_j(x_j)|^2+
|g_j(( I_Ux)_j(x_j))|^2\}\,d\mu_j(x_j)\,dt=I_{j,n}^++I_{j,n}^-,
\]
where
\beqs
I_{j,n}^+&:=&
\int_{\Sigma_{j,n}^+}\{|( I_Ux)_j(x_j)|^2+
|g_j(( I_Ux)_j(x_j))|^2\}\,d\mu_j(x_j)\,dt,
\\
I_{j,n}^-&:=&
\int_{\Sigma_{j,n}^-}\{|( I_Ux)_j(x_j)|^2+\alpha_j(t,x_j)^2
|g_j(( I_Ux)_j(x_j))|^2\}\,d\mu_j(x_j)\,dt.
\eeqs

For the estimation of $I_{j,n}^+$, we first  notice that
the assumption  (\ref{gbounded})  leads to
\[
I_{j,n}^+\leq (1+2M)
\int_{\Sigma_{j,n}^+}( |( I_Ux)_j(x_j)|^2\,d\mu_j(x_j)\,dt,
\]
and  by (\ref{gstrp})  and \eqref{upperbdalphaj} we get
\[
I_{j,n}^+\leq   m^{-1} (1+ 2M) \int_{\Sigma_{j,n}^+}  ( I_Ux)_j(x_j)\cdot g_j(( I_Ux)_j(x_j)) \,d\mu_j(x_j)dt. 
\]
By  the assumption \eqref{alphaoscj}, we directly obtain 
\[
I_{j,n}^+\leq   m^{-1} (1+ 2M)\alpha_0^{-1} \int_{\Sigma_{j,n}^+}   \alpha_j(t, x_j) ( I_Ux)_j(x_j)\cdot g_j(( I_Ux)_j(x_j)) \,d\mu_j(x_j)dt.
\]
By \eqref{serge:11/10:1},
we arrive at
\be\label{spo10osc}
I_{j,n}^+
\leq m^{-1} (1+ 2M)\alpha_0^{-1}   (\caE(nT)-\caE((n+1)T)).
\ee

Similarly by the assumption (\ref{gnear0})
and the monotonicity of $G$ and $\alpha$ we have 
\[
I_{j,n}^- \leq   \int_{nT}^{(n+1)T}
\int_{X_j} G(( I_Ux)_j(x_j)\cdot g_j(( I_Ux)_j(x_j))) \,d\mu_j(x_j)dt.
\]
Jensen's inequality then yields
\[
I_{j,n}^-\leq    T \mu_j(X_j) 
G\left(\frac{1}{T\mu_j(X_j)}  \int_{nT}^{(n+1)T}
\int_{X_j}  ( I_Ux)_j(x_j)\cdot g_j(( I_Ux)_j(x_j)) \,d\mu_j(x_j)dt
\right).
\]
As  $G$ is strictly increasing and again using \eqref{alphaoscj}, we obtain 
 \beqs
I_{j,n}^-\leq   
T \mu_j(X_j) 
G\left(\frac{1}{T\mu_j(X_j)\alpha_0} \!\! \int_{nT}^{(n+1)T}
\int_{X_j} \alpha_j(t,x_j)( I_Ux)_j(x_j)\cdot g_j(( I_Ux)_j(x_j)) \,d\mu_j(x_j)dt
\right).
\eeqs
By  (\ref{2.2weak}), we arrive at
\be\label{spo11osc}
I_{j,n}^-\leq \ T \mu_j(X_j) 
G\left(\frac{\caE(nT)-\caE((n+1)T)}{T \mu_j(X_j)\alpha_0}
\right).
\ee

The estimates   (\ref{spo10osc}) and  (\ref{spo11osc})  into the estimate
(\ref{spo9bisalphadecroissant}) and the monotonicity of $G$ give
 \[
\caE((n+1)T)
\leq c_2    \left\{\frac{\caE(nT)-\caE((n+1)T)}{T\mu \alpha_0}
+G\left(\frac{\caE(nT)-\caE((n+1)T)}{T\mu \alpha_0}\right)\right\},
\]
for some positive constant $c_2$ (depending on $T$, $C$, $\omega$ (from  Theorem \ref{tec}),  $\max_j\mu_j(X_j)$, $\tilde \alpha$, $\alpha_0$, $M$ and $m$).
Hence,
\beqs
\caE(nT)&=&\caE(nT)-\caE((n+1)T)+\caE((n+1)T)
\\
&\leq& \max\{T\mu \alpha_0,c_{2}\}
\left\{\frac{\caE(nT)-\caE((n+1)T)}{T\mu \alpha_0}
+G\left(\frac{\caE(nT)-\caE((n+1)T)}{T\mu \alpha_0}\right)\right\}.
\eeqs
With $c=\max\{T\mu \alpha_0,c_{2}\}$, and the definition \eqref{defh} of $h$, we have found that
%
%
%
%
\[
T\mu  \alpha_0 h^{-1}(\caE(nT))+\caE((n+1)T) \leq \caE(nT).
\]
Since this estimate is valid for all  $n \in \mathbb{N}$, we conclude by Theorem \ref{tii}
with the choice $\beta(t)=T\mu  \alpha_0$.
\end{proof}

\bex\label{exempleosc}
{\rm
If $g_j$ satisfies the assumptions from Example \ref{exemple}, then in the linear case (i.e., if $p=\gamma=1$) 
we get the exponential decay
 decay 
 \[
 \caE(t)\leq K  \caE(0) e^{-Lt}, \forall t\geq 0,
 \]
 for some positive constants $K$ and $L$.  
On the contrary if $p+1>2\gamma$, then  we get  the decay
$K(\caE(0)) t^{-\frac{2\gamma}{p+1-2\gamma}}$. In both cases, the decay rate is the same as the one of the autonomous case.
}
\eex
 
\section{Illustrative examples\label{sex}}

 \subsection{Second order evolution equations}\label{s2ordereveq}

Some examples given below enter in the following framework:
Let $H$ and $V$ be  two real separable  Hilbert spaces    such that
$V$ is densely and continuously embedded into $H$.
Define the linear operator $A_2$ from $V$ into $V'$ by
\be\label{a2ip}
\langle A_2 u,v\rangle_{V'-V}=(u,v)_V, \forall u,v \in V,
\ee
and suppose given  a (nonlinear) and time-dependent mapping $B_2(t)$ from $V$ into $V'$
as follows: We assume that $V$ is continuously embedded into   a  control space $U$
in the form \eqref{defU} with the same assumptions on $U_j$, $j=1,\cdots, J$.  
Similarly,  we suppose given    mappings $g_j$ 
and $\alpha_j$ satisfying the same assumptions than in subsection \ref{ssgeneralhypo}.
We then  define  the (nonlinear)  operator $B_2(t)$ from $V$ into $V'$   by
\be\label{defB2nonlineaire}
\langle B_2(t) u,v\rangle=  \sum_{j=1}^J\int_{X_j} \alpha_j(t,x)  g_j((J_Uu)_j(x_j))\cdot ( J_Uv)_j(x_j)\,d\mu_j(x_j), \forall u,v\in V,
\ee
where  $J_U$ denotes the embedding from $V$ to $U$ 
(hence, $(J_Uu)_j$ is the $j^{th}$ component of $J_Uu$).


With these data, we
consider   the second order evolution equation
\be\label{soeveq}
\left\{
\begin{tabular}{lllll}
& $\frac{d^2 u}{dt^2}(t)+A_2 u(t)+ B_2(t)\frac{du}{d t}(t)=0 $ in $H, t\geq 0$,\vspace{2mm}
\\
& $u(0) =  u_0, \frac{du}{d t}(0)=u_1.$
\end{tabular}
\right.
\ee

This system is reduced to the first order system (\ref{a6})
using the standard argument of reduction of order: setting
$\caH=V\times H$, $\caV=V\times V$
with their natural inner products,
$$
x=(u,v)^\top,
$$
with $v=\frac{du}{d t}$
and introducing the operators
\be\label{defA1andB(t)roo}
A_1 x=(-v,A_2u)^\top, B(t)x= (0,B_2(t)v)^\top.
\ee
Note that $B(t)$ is indeed in the form \eqref{defBnonlineaire} with 
  $I_U(u,v)^\top=J_U v$, for all $(u,v)^\top \in V\times V$.

With this definition, we see that $x$ is solution of \eqref{a6},  assuming that $u$ exists and is sufficiently regular. 
But in its full generality, the domain of $A_B(t)$ is time-dependent; so, again we distinguish between two cases.

Before going on, let us notice that the above operator $A_1$ trivially satisfies 
\eqref{h4} due to \eqref{a2ip}. Consequently the (nonlinear) operator $A_B(t)=A_1 x+B(t)$
 corresponding to \eqref{defA1andB(t)roo} satisfies all assumptions of subsection \ref{ssgeneralhypo}.

Let us finally remark that Theorem 6.1 of 
\cite{nicaise:rendiconti:03}   shows that $A_1-I_U$ and $-A_1-I_U$ generates a $C_0$-semigroup of contractions in $\caH$.

\subsubsection{The  bounded case}

\bt\label{rkB_2(t)}
In addition to the above assumptions, if we  assume that $H$ is continuously embedded into   the  control space $U$, then for all $(u_0, u_1)\in   D(A_2)\times V$  problem \eqref{soeveq} has a unique solution $u\in
C([0,\infty), V)\cap C^1([0,\infty), H)$ such that its second derivative $u''(t)=-A_2 u(t)-B_2(t)u'(t)$ exists  and is continuous in H, except at a countable number of points $t$.
\et
\begin{proof}
We show that  $A_B(t)=A_1 x+B(t)$ satisfies the assumptions of Theorem \ref{texistencebdB}. First as $\caH=V\times H$, it is clearly embedded into $U$ as $H\hookrightarrow U$ and that $D(A_B(t))=D(A_2)\times V$. Hence, it suffices to show that there exists a positive real number $\lambda$ such that  $R(\lambda \mathbb{I}+\caA_B(t))=\caH$. But this properties is proved in \cite[Theorem 6.1]{nicaise:rendiconti:03} for $\lambda=1$.
We then conclude by Theorem \ref{texistencebdB} that for any $(u_0, u_1)\in D(A_2)\times V$, there exists 
a unique solution 
$x$ of \eqref{a6} with the properties  \eqref{regsolutiona6}. 
We now come back to the original system by noticing that $x(t)=(u(t), v(t))^\top$ satisfies
\[
(u'(t), v'(t))^\top=(v(t), -A_2u(t)-B_2(t)v(t))^\top, \forall t\geq 0.
\]
Hence, $u\in C^1([0,\infty), H)$, $v=u'$ and the second components of the above identity yields $u''(t)=-A_2 u(t)-B_2(t)u'(t)$.

The proof is complete.
\end{proof}

 \subsubsection{The unbounded case}

 Here in order to avoid the time-dependency of the domain of $A_B(t)$, we suppose that
 the mappings $\alpha_j$ satisfies \eqref{alphaequal} for some $\alpha\in C([0,\infty), (0,\infty))$. In such a case, the operator 
 $B_2(t)$ defined in \eqref{defA1andB(t)roo} will be in the form
 $B(t)=\alpha(t) B_1$, where $B_1(u,v)^\top=(0,B_2v)^\top$, with  (compare with 
 \eqref{defBnonlineairepartcase})
\be\label{defB2nonlineairetimeind}
\langle B_2 u,v\rangle=  \sum_{j=1}^J\int_{X_j}   g_j((J_Uu)_j(x_j))\cdot ( J_Uv)_j(x_j)\,d\mu_j(x_j), \forall u,v\in V.
\ee

Under the previous assumptions on  $A_2$  and $B_2$, with the help of Theorem \ref{texistencebubdB} we can prove the next existence result for problem  \eqref{soeveq}.
\bt\label{tsoeveq}
In addition to the above assumptions, we assume that the mappings $\alpha_j$ satisfies \eqref{alphaequal} for some $\alpha\in C^1([0,\infty), (0,\infty)$ such that $\alpha'$ is locally Lipschitz. Then
for all $(u_0, u_1)\in V\times V$ such that
$A_2 u_0+\alpha(0)B_2u_1\in H$, problem \eqref{soeveq} has a unique solution $u\in
C([0,\infty), V)\cap C^1([0,\infty), H)$ such that its second derivative $u''(t)=-A_2 u(t)-\alpha(t)B_2u'(t)$ exists  and is continuous in H, except at a countable number of points $t$.
\et
\begin{proof}
We first recall that $
x=(u,v)^\top
$ is solution of
 (\ref{a6}) with $A_1$ and $ B(t)$ from \eqref{defA1andB(t)roo} (and $B_2(t)=\alpha(t) B_2$) if and only if
 \beqs
 u'(t)&=& v(t),\\
 v'(t)&=&-A_2 u(t)-\alpha(t) v(t).
 \eeqs
We now perform the following change of unknowns (assuming that $u, v$ exists and are sufficiently regular)
\be\label{cgtvariables}
\tilde u(t)= \alpha(t)^{-1} u(t), \quad \tilde v(t)=v(t).
\ee
Then setting 
$\tilde x=(\tilde u(t), \tilde v(t))^\top$, we see that 
it satisfies
\be\label{equationx'=Ax}
\tilde x'=\left(
\begin{array}{ll}
-\alpha(t)^{-2}\alpha'(t) u(t)+\alpha(t)^{-1} v(t)\\
-(A_2 u(t)+\alpha(t)B_2 v(t))
\end{array}
\right)
=\left(
\begin{array}{ll}
\alpha(t)^{-1}(-\alpha'(t) \tilde u(t)+\tilde v(t)\\
- \alpha(t) (A_2 \tilde u(t)+B_2 \tilde v(t))
\end{array}
\right).
\ee
This means that as operator $D(t)\in \caL (\caH)$, we here choose 
\be\label{defD}
D(t) (u,v)^\top= (\alpha(t)^{-1} u,   v)^\top, \forall (u,v)^\top\in V\times H.
\ee
From the previous identity \eqref{equationx'=Ax}, we see that the assumption 
\eqref{commutateur} holds with
\be\label{deftildeD}
\tilde D(t) (u,v)^\top= (u,  \alpha(t) v)^\top, \forall (u,v)^\top\in V\times H.
\ee
From the assumptions on $\alpha$ and their definitions, we readily check that all other assumptions from Theorem \ref{texistencebubdB} on $D$ and $\tilde D$ are satisfied.
Finally   Theorem 6.1 of  \cite{nicaise:rendiconti:03} (since $B_2$ defined above satisfies the assumption of this Theorem) guarantees that $A_1+B_1$ is maximal monotone
and has a dense domain in $\caH$.
 In conclusion, by Theorem \ref{texistencebubdB}, if 
 $(u_0, u_1)\in D(A_B(0))$ (or equivalently if $(u_0, u_1)\in V\times V$ is such that
$A_2 u_0+\alpha(0)B_2u_1\in H$),   there exists a unique solution 
$x$ of \eqref{a6} with the properties  \eqref{regsolutiona6}.
\end{proof}

In the remainder of this section
 $\Omega $ is a bounded domain of $\R^n$, $n\geq 1$ with a Lipschitz boundary  $\Gamma $.
Some restrictions will be specified later on  when they will be necessary.
We further denote by $\nu$ the unit outward normal vector along $\Gamma$.

\subsection{Nonlinear and nonautonomous stabilization of the wave equation}\label{subsectwave}

\subsubsection{Interior damping\label{sswaveint}}

Consider the wave equation with interior damping and Dirichlet boundary condition
\be\label{wave}
\left\{
\begin{tabular}{llllllll}
& $\partial_t^2 u-\Delta u+ \sigma  \sum_{j=1}^J \alpha_j(t,\cdot) g_j(\partial_t u)=0 $ in $Q:=\Omega \times ]0,+\infty[$,\vspace{2mm}
\\
&$u=0$ on $\Sigma:= \Gamma\times ]0,+\infty[$,\vspace{2mm}
\\
& $u(0) =  u_0,$
$\partial_t u(0) =  u_1$ in $\Omega $,\vspace{2mm}
\end{tabular}
\right.
\ee
where  $\sigma$ is a non-negative function that belongs to $L^\infty(\Omega)$
such that 
that there exists a positive constant $\sigma_0$ such that
\be \label{deltacond}
\sigma\geq \sigma_0 \hbox{ on }\mathcal{O},
\ee
for some open and non empty subset $\mathcal{O}$ of the the support $X_\sigma$
 of $\sigma$. For all $j=1,\cdots, J$, the functions $\alpha_j$ and $g_j$ satisfy the assumptions of subsection \ref{ssgeneralhypo} with $U_j=L^2(X_j)$, $X_j$ being an open 
 and non empty subset of $X_\sigma$ such that 
 \be\label{partition}
 X_j\cap X_k=\emptyset, \hbox{ for } j\ne k, \hbox{ and } \cup_{j=1}^J \bar X_j=X_\sigma.
 \ee

The stability of this problem in the autonomous case, namely for $\alpha_j=1$, was extensively studied in the litterature,
let us cite the papers  \cite{haraux:89a,Komornik93a,lebeau:96,Liuzuazua,martinez:99,zuazua:cpde:90} 
and the references cited there. Both papers are restricted to some particular choices
of $\sigma$   and $g_j$ leading to some exponential or polynomial decay rates of
the energy of the solution of (\ref{wave}).
On the contrary the nonautonomous case is less considered in the literature
and with the exception of \cite{mustapha:15} all papers concern interior damping acting on the whole domain (i.e. $\sigma=1$), see
\cite{Daoutli:11,Jiao-Xiao:15,Luo-Xiao:20a,Luo-Xiao:20b,Luo-Xiao:21,martinez:00,Nakao:97,Pucci-Serrin:96}.
Using the results of the previous section,
and under the assumption that the autonomous  linear system is exponentially stable,
we  obtain new decay results  for a large class of functions  $g_j$   
and $\alpha_j$.

The first point is that problem (\ref{wave})
 enters in the framework of
problem  \eqref{soeveq}
from subsection \ref{s2ordereveq} once we take
\beqs
H&=&L^2(\Omega),\\
V&=&H^1_0(\Omega),\\
(u,v)_V&=&\int_\Omega \nabla u\cdot \nabla v\, dx,  \forall u,v\in V,\\
\langle B_2(t) u,v\rangle_{V'-V}&=& \sum_{j=1}^J \int_{X_j} \alpha_j(t,x) \sigma(x) f_j(u(x)) v(x)\,dx, \forall u,v\in V.
\eeqs

Let us notice that the inner product
$(\cdot,\cdot)_V$ induces a norm on $V$
equivalent to the usual one due to Poincar\'e inequality.  Furthermore, the condition \eqref{gbounded} allows to show that $B_2(t)$ is well-defined from $V$ to $V'$.

As $L^2(\Omega)$ is clearly embedded into $U=\prod_{j=1}^J L^2(X_j)$ (that is clearly identical with $L^2(X_\sigma)$), the assumptions of Theorem  \ref{rkB_2(t)} are satisfied and therefore, there exists a unique solution 
$u$ of (\ref{wave})  such that $(u,u')^\top$ satisfies \eqref{regsolutiona6}.

In order to deduce some stability results for our system (\ref{wave}) with the help of Theorem \ref{tstabobsgen}
 we need  that  $-A_1-I_U$ generates an exponentially stable semigroup  in $\caH$,
 with
 the control space $U=L^2(X_\sigma)$. This property is equivalent to the exponential decay
of the solution  of the autonous and linear problem
 \be\label{waveint}
\left\{
\begin{tabular}{llllllll}
& $\partial_t^2 u-\Delta u+  \sigma   \partial_t u=0 $ in $Q:=\Omega \times ]0,+\infty[$,\vspace{2mm}
\\
&$u=0$ on $\Sigma$,\vspace{2mm}
\\
& $u(0) =  u_0,$
$\partial_t u(0) =  u_1$ in $\Omega $.\vspace{2mm}
\end{tabular}
\right.
\ee

 Note that the exponential stability of \eqref{waveint} holds in many different situations, see \cite{haraux:89a,zuazua:cpde:90} in the case of a $C^2$ boundary 
and $\mathcal{O}$ being a neighborhood of 
\be\label{defG+}
\Gamma_+:=\{x\in \Gamma: (x-x_0)\cdot \nu(x)>0\},
\ee
for some $x_0\in \R^n$,
 or
  \cite{lebeau:96} in the case of a domain $\Omega$ with an analytical boundary,  $\sigma$ smooth
and $\mathcal{O}$ satisfying a geometrical control condition.
Note that in the case $d=1$, this assumption is valid as soon as
$\mathcal{O}$ contains an open interval of $\Omega$, see  \cite[Exemple 1]{haraux:89a}.
Moreover, if the linear damping acts on the whole domain, namely if $\sigma= 1$ in 
 \eqref{waveint} a simple spectral analysis shows that \eqref{waveint} is exponentially stable without any assumption on the regularity of the boundary of $\Omega$.
In all these situations, if $g_j$ and
$\alpha_j$ satisfy the additional assumptions of Theorem \ref{tstabobsgen}, \ref{tstabobsgencroissant} or \ref{tstabobsgenosc}, then 
the energy of our system will satisfy \eqref{A19}, \eqref{A19croissante} or \eqref{A19osc}.
This allows to recover and extend some results from \cite{Pucci-Serrin:96,Nakao:97,martinez:00,Jiao-Xiao:15,Daoutli:11,Luo-Xiao:20a,Luo-Xiao:20b,Luo-Xiao:21}. Particular cases not covered by the previous references are the case when we have only a local damping, namely  $X_\sigma\ne \bar \Omega$, and/or a factor $\alpha(t)$ piecewise variables, for instance
\[
\alpha_j(t,x)=\alpha_j(t),\forall x\in X_j.
\]

\subsubsection{Boundary damping\label{ss522}}

Consider the wave equation with a boundary damping
\be\label{wavebdy}
\left\{
\begin{tabular}{llllllll}
& $\partial_t^2 u-\Delta u=0 $ in $Q:=\Omega \times ]0,+\infty[$,\vspace{2mm}
\\
&$u=0$ on $\Sigma_0:= \Gamma_0\times ]0,+\infty[$,\vspace{2mm}
\\
&$\partial_\nu u+ au + \alpha(t)  k(x) g(\partial_t u)=0$ on $\Sigma_1:= \Gamma_1\times ]0,+\infty[$,\vspace{2mm}
\\
& $u(0) =  u_0,$
$\partial_t u(0) =  u_1$ in $\Omega $,\vspace{2mm}
\end{tabular}
\right.
\ee
where $\Gamma_0$ is an open subset of $\Gamma$, 
 $\Gamma_1=\Gamma\setminus \bar\Gamma_0$,  $a, k\in L^\infty(\Gamma_1)$ are two non negative real-valued functions.  The function  $g$ is  a non-decreasing continuous function
from $\R$ into itself such that $g(0)=0$ and satisfying \eqref{gbounded}, while the function   $\alpha\in C^1([0,\infty; (0,\infty))$ 
and $\alpha'$ is   locally Lipschitz.

 For the sake of simplicity we suppose that
\be\label{assgamma0a}
\hbox{ either } \Gamma_0 \hbox{ is not empty or } a\not\equiv 0.
\ee

As previously, the stability of this problem in the autonomous case, namely for $\alpha=1$, was extensively studied in the litterature,
let us cite the papers  
\cite{bardos,Chen_G:79b,Chen_G:79a,Komornik91,Komornik93b,Komornikbook,KomornikZuazua,LaJDE83,Lasiecka-Triggiani92,Triggiani:89,zuazua:sicon:90}
and the references cited there. Both papers are restricted to some particular choices
of $\Gamma_0$, $a$, and $g$ leading to some exponential or polynomial decay rates of
the energy of the solution of (\ref{wave}).
On the other hand to the best of our knowledge  the nonautonomous case is only considered in  \cite{mustapha:15}.
Using the results of the previous section,
and under the assumption that the autonomous  linear system is exponentially stable,
we  obtain new decay results  for a large class of functions  $g$ 
and $\alpha$.

As before problem (\ref{wavebdy})
 enters in the framework of
problem (\ref{soeveq})
from subsection \ref{s2ordereveq} once we take:
\beqs
H&=&L^2(\Omega),\\
V&=&\{v\in H^1(\Omega)| v=0 \hbox{ on } \Gamma_0\},
\\
(u,v)_V&=&\int_\Omega \nabla u\cdot \nabla v\, dx+\int_{\Gamma_1}
 a u\cdot v\,d\sigma, \forall u,v\in V,\\
U&=&L^2(\Gamma_1), \\
\langle B_2(t) u,v\rangle_{V'-V}&=&\alpha(t) \int_{\Gamma_1} k(x) g(u(x)) v(x)\,d\sigma(x), \forall u,v\in V.
\eeqs

Let us remark that the assumption (\ref{assgamma0a}) implies that the inner product
$(\cdot,\cdot)_V$ induces a norm on $V$
equivalent to the usual one, while our condition    \eqref{gbounded}   implies that  $B_2(t)$ is well-defined.

We readily check that these assumptions guarantee that $B_2(t)$ fulfils
all the assumptions of Theorem \ref{tsoeveq};  hence, (\ref{wavebdy})  has  a unique solution 
$u$   such that $(u,u')^\top$ satisfies \eqref{regsolutiona6}.

In order to take advantage of Theorem \ref{tstabobsgen}
 we need  that  $-A_1-I_U$ generates an exponentially stable semigroup   in $\caH$. For this particular example this property is equivalent to the exponential decay
of the solution  of the autonous and linear problem
 \be\label{wavbdy}
\left\{
\begin{tabular}{llllllll}
& $\partial_t^2 u-\Delta u=0 $ in $Q:=\Omega \times ]0,+\infty[$,\vspace{2mm}
\\
&$u=0$ on $\Sigma_0:= \Gamma_0\times ]0,+\infty[$,\vspace{2mm}
\\
&$\partial_\nu u+ au + k\partial_t u=0$ on $\Sigma_1:= \Gamma_1\times ]0,+\infty[$,\vspace{2mm}
\\
& $u(0) =  u_0,$
$\partial_t u(0) =  u_1$ in $\Omega $.\vspace{2mm}
\end{tabular}
\right.
\ee

The   exponential stability of \eqref{wavbdy} was obtained in many different situations, 
let us quote  
\cite{Chen_G:79a,Chen_G:79b}, where $a=0$, $k\in L^\infty(\Gamma_1)$
such that
\be\label{conk}
k\geq k_0 \hbox{ on } \Gamma_1,
\ee
for some positive constant $k_0$
and under the assumptions that
\beq
 &&m\cdot \nu\leq 0 \hbox{ on } \Gamma_0,
\label{condchen0}
\\
\label{condchen1}
&& m\cdot \nu\geq \gamma>0 \hbox{ on } \Gamma_1,
\eeq
where $\gamma$ is a positive constant and    $m$ is the standard multiplier defined by
\be\label{multstandard}
m(x)=x-x_0, \forall x\in \R^n,
\ee
for some point $x_0\in \R^n$. 

This result was generalized in \cite{LaJDE83,Triggiani:89} to a more general class of 
multipliers $m\in C^2(\bar \Omega)$ for which
the matrix $(\partial_j m_i+\partial_i m_j)_{1\leq i,j\leq n}$ is uniformly positive definite in 
 $\bar \Omega$ but still under the assumptions $a=0$, $k\in L^\infty(\Gamma_1)$
satisfying \eqref{conk} and the geometrical constraints \eqref{condchen0}-\eqref{condchen1}.

Let us observe that conditions  \eqref{condchen0}-\eqref{condchen1} force to have
\be\label{emptyintersection}
\bar \Gamma_0\cap \bar \Gamma_1=\emptyset.
\ee
This constraint has been removed in
\cite{Lasiecka-Triggiani92} since condition \eqref{condchen1} has been removed, while the other conditions from \cite{LaJDE83,Triggiani:89}  remain. Alternatively, in  \cite{KomornikZuazua,Lagnese:88}, the choice 
$k=m \nu$ (with $m$ in the form \eqref{multstandard} and then as in \cite{LaJDE83,Triggiani:89}) allows to replace the condition 
\eqref{condchen1} by
\[
 m\cdot \nu> 0 \hbox{ on } \Gamma_1,
\]
under the conditions  $a=0$ and $\Gamma_0$  non empty, see also \cite{Komornik91,martinez:99} for the case $a\not\equiv 0$. 

Let us finally notice that microlocal analysis arguments from \cite{bardos}  allow to suppress the condition 
\eqref{condchen0} if $\Gamma$ is analytic, the condition \eqref{emptyintersection} holds, $a$ and $k$ are smooth, and $\Gamma_1$ satisfies the geometrical control condition that it must meet each ray in a nondiffractive point.

Since in 
\cite{mustapha:15}, it is assumed that $a=0$, $k=1$, that  \eqref{condchen0}-\eqref{condchen1} hold with $m$ in the form \eqref{multstandard} and that \eqref{emptyintersection} holds, Theorems \ref{tstabobsgen}, \ref{tstabobsgencroissant} and \ref{tstabobsgenosc} allow to improve significantly the result from \cite{mustapha:15} by obtaining different decay rates of the solution of system \eqref{wavebdy} with appropriated choices of $\alpha$ and $g$
using the above mentioned results about the exponential decay of system \eqref{wavbdy}.

 \subsubsection{Pointwise interior damping}
 
In this subsection, we conisder  the large time behavior of the solution of  
a homogenous string equation with a
 homogenous Dirichlet boundary condition at the left end 
and a Neuman boundary condition at the right end subject to 
a time-dependent and nonlinear pointwise interior actuator. More precisely, we condider the problem
\be \label{pointwise} 
\left\{
\begin{tabular}{llllllll}
& $\partial^2_t u  - \, \partial^2_x u  +  
\alpha(t) g(\partial_t u)  \, \delta_{\xi} = 0$  in $(0,\pi)\times \R,$
\\
&$u(0,t) =  \partial_x u(\pi,t) = 0, \, t > 0,$
\\
&$u(\cdot,0) = u_0, \, \partial_t u(\cdot,0) = u_1$  in $ (0,\pi), $
\end{tabular}
\right.
\ee
where $\xi$ is a fixed point of $(0,\pi)$, the functions  $g$ is  a non-decreasing continuous function
from $\R$ into itself such that $g(0)=0$  and satisfying \eqref{gbounded}, and the function   $\alpha\in C^1([0,\infty; (0,\infty))$ is such that
  $\alpha'$ is   locally Lipschitz.
 
 The stability of this problem in the autonomous and linear case, namely for $\alpha=g=1$ was considered in
\cite{ammari:01b} (see also \cite{ammari:01}), and to the best of our knowledge, the case of a nonautonomous and nonlinear  pointwise damping has not been analyzed.

Let us notice that problem (\ref{pointwise})
 enters in the framework of
problem (\ref{soeveq})
from subsection \ref{s2ordereveq} once we take:
\beqs
H&=&L^2(0,\pi),\\
V&=&\{v\in H^1(0,\pi)| v(0)=0\},
\\
(u,v)_V&=&\int_0^\pi   u_x  v_x\, dx, \forall u,v\in V,\\
U&=&\R, \\
\langle B_2(t) u,v\rangle_{V'-V}&=&\alpha(t)  g(u(\xi)) v(\xi),  \forall u,v\in V.
\eeqs

These assumptions guarantee that $B_2(t)$ fulfils
all the assumptions of Theorem \ref{tsoeveq};  hence, (\ref{pointwise})  has  a unique solution 
$u$   such that $(u,u')^\top$ satisfies \eqref{regsolutiona6}.

As Theorem 1.2 of \cite{ammari:01b} guarantees the exponential decay of the solution of 
\eqref{pointwise} with $\alpha=g=1$ if $\frac{\xi}{\pi}=\frac{p}{q}$ with $p\in \N^*$ odd and $q\in \N^*$, we can apply Theorem   \ref{tstabobsgen}, \ref{tstabobsgencroissant} or \ref{tstabobsgenosc} to  obtain different decay rates of the solution of system \eqref{pointwise} under this assumption on $\xi$ and  if 
$\alpha$   and $g$ satisfy the additional assumptions  from Theorem \ref{tstabobsgen}, \ref{tstabobsgencroissant} or \ref{tstabobsgenosc}.

\subsection{Nonlinear and nonautonomous stabilization of the elastodynamic system}

%
%
%
%
%

With the notation of the above subsubsection \ref{ss522}, we consider the following elastodynamic system:
\be\label{elasticite}
\left\{
\begin{tabular}{llllllll}
& $\partial_t^2 u-\nabla \sigma (u)+ \sigma  \sum_{j=1}^J \alpha_j(t,\cdot) g_j(\partial_t u)=0 $ in $Q$,\vspace{2mm}
\\
& $u=0 $
on $\Sigma_0$,\vspace{2mm}
\\
& $\sigma (u)\cdot \nu+a u+kg(\partial_t u)=0 $
on $\Sigma_1$,\vspace{2mm}
\\
& $u(0) =  u_0,$
$\partial_t u(0) =  u_1$ in $\Omega $.
\end{tabular}
\right.
\ee
As usual $u(x,t)$ is the displacement field at the point $x\in\Om$ at time $t$
and  $\sigma (u)=(\sigma_{ij} (u))_{i,j=1}^n$ is the stress tensor given by (here  and in the sequel we shall use the summation convention for repeated indices)
$$
\sigma_{ij} (u)=a_{ijkl}\eps_{kl}(u),
$$
where $\eps (u)=(\eps_{ij} (u))_{i,j=1}^n$ is the strain tensor given by
$$
\eps_{ij} (u)=\frac12(\frac{\partial u_i}{\partial x_j} +\frac{\partial u_j}{\partial x_i}),
$$
and the tensor $(a_{ijkl})_{i,j,k,l=1,\cdots, n}$ is made of $W^{1,\infty}(\Om)$ entries
such that
$$
a_{ijkl}=a_{jikl}=a_{klij},
$$
and
 satisfying the ellipticity condition
$$
a_{ijkl} \eps_{ij}\eps_{kl}\geq \alpha  \eps_{ij}\eps_{ij},
$$
for every symmetric tensor $( \eps_{ij})$ and some $ \alpha>0$.
Hereabove and below $\nabla \sigma (u)$ is the vector field defined by
$$
\nabla \sigma (u)=(\partial_j \sigma_{ij} (u))_{i=1}^n.
$$

Finally $a$ and $k$ are two nonnegative real number.  As before we   assume that
 \be\label{hopFG}
g_j=0, \forall j=1,\cdots, J \hbox{ or }  \Gamma_1=\emptyset.
\ee
This last assumption means that we stabilizate our system either by a boundary feedback
or by an internal feedback with only Dirichlet boundary conditions.
In case of a boundary damping,   we also  suppose  that (\ref{assgamma0a})
 holds.

In case of an internal feedback,
as in subsubsection \ref{sswaveint},  $\sigma$ is a non-negative function that belongs to $L^\infty(\Omega)$ satisfying \eqref{deltacond} 
for some open and non empty subset $\mathcal{O}$ of the the support $X_\sigma$
 of $\sigma$. For all $j=1,\cdots, J$, the functions $\alpha_j$ and $g_j$ satisfy the assumptions of subsection \ref{ssgeneralhypo} with $U_j=L^2(X_j)^n$, $X_j$ being an open  and non empty
 subset of $X_\sigma$ such that \eqref{partition} holds.

In case of a boundary feedback, the functions $\alpha$ and $g$ satisfy the assumptions of subsection \ref{ssgeneralhypo} with $U=L^2(\Gamma_1)^n$
and suppose, moreover, that  $\alpha\in C^1([0,\infty; (0,\infty))$ is such that
  $\alpha'$ is   locally Lipschitz.

The stability of the  system (\ref{elasticite}) was considered in
\cite{AlabauKomornik,Bey:03,guesmia98b,Guesmia99b,guesmia:00,Horn,LaSiam83,Tebou:96}
in the autonomous case
under  some particular hypotheses
on $\Gamma_0$,  $\Gamma_1$, $a$, $g_j$ and $g$ leading to  exponential or polynomial decay   of
the energy of the solution of (\ref{elasticite}). The nonautonomous case  with internal feedback and for the Lam\'e system (corresponding to $n=3$ and to the choice $a_{ijkl}=\lambda \delta_{ij}\delta_{kl}+\mu (\delta_{ik}\delta_{jl}+\delta_{il}\delta_{jk})$, where $\lambda$ and $\mu$ are  positive constants, called Lam\'e parameters)  was treated in
\cite{BchatniaDaoulatli:13,Bellassoued:08}.

As in the above subsection,  problem (\ref{elasticite})
 may be expressed in the form  (\ref{soeveq})
from subsection \ref{s2ordereveq} with the choices:
\beqs
H&=&L^2(\Omega)^n,\\
V&=&\{v\in H^1(\Omega)^n| v=0 \hbox{ on } \Gamma_0\},
\\
(u,v)_V&=&\int_\Omega \sigma_{ij}(u)\epsilon_{ij}(v)\,
dx+a\int_{\Gamma_1} u\cdot v\,d\sigma, \forall u,v\in V,
\eeqs
and
\[
\langle B_2(t) u,v\rangle_{V'-V}
=\sum_{j=1}^J \int_{X_j} \alpha_j(t,x) \sigma(x) f_j(u(x))\cdot v(x)\,dx, \forall u,v\in V,
\]
in case of an interior damping and
\[
\langle B_2(t) u,v\rangle_{V'-V}=\alpha(t)\int_{\Gamma_1} g(u)\cdot v\,d\sigma, \forall u,v\in V,
\]
otherwise.

In the case of an interior damping (resp. boundary  damping), all the assumptions of Theorem  \ref{rkB_2(t)} (resp. Theorem  \ref{tsoeveq}) are satisfied and therefore, we  have a unique solution 
$u$ of (\ref{elasticite})  such that $(u,u')^\top$ satisfies \eqref{regsolutiona6}.

For stability results of \eqref{elasticite},
 we need to check that  $-A_1-I_U$ generates an exponentially stable semigroup in $V\times H$,
where the control space $U$ is  defined by
\beqs
U&=&L^2(X_\sigma)^n \hbox{ if } \Gamma_1=\emptyset,\\
U&=&L^2(\Gamma_1)^n   
\hbox{ if }  g_j=0, \forall j=1,\cdots, J.
\eeqs
As before, this is equivalent to the exponential decay of the autonomous and linear system  (\ref{elasticite}), i.e. corresponding to
   $\Gamma_1=\emptyset$, $\alpha_j=1$ and $g_j(s)=s$ in the first case and
to $g_j=0$, $\alpha=1$ and $g(s)=s$ in the second case.

In the first case (i. e.,   $\Gamma_1=\emptyset$), this exponential decay   was proved in  \cite[Theorem 1.1]{guesmia:00} (see also \cite{guesmia98b} for the case $X_\sigma=\Omega$) under the assumptions that 
$\mathcal{O}$ is a neighborhood of $\Gamma_+$ defined by
\eqref{defG+}. Hence, in the setting of one of these papers, under the additional assumptions    on $\alpha_j$ and $g_j$ from Theorem \ref{tstabobsgen}, \ref{tstabobsgencroissant} or \ref{tstabobsgenosc},   different decay rates of the solution of system \eqref{elasticite} (with $\Gamma_1=\emptyset$) are available.

In the second case (i.e., $g_j=0,$ for all $j=1,\cdots, J$), the exponential decay of the autonomous and linear system  (\ref{elasticite})   was proved  in
\cite{AlabauKomornik,Bey:03,Horn,LaSiam83} under some geometric assumptions. In the setting of one of these papers, we then obtain different decay rates of the solution of system \eqref{elasticite} (with $g_j=0,$ for all $j=1,\cdots, J$) if   $g$ 
 and $\alpha$ satisfy the assumptions   from Theorem \ref{tstabobsgen}, \ref{tstabobsgencroissant} or \ref{tstabobsgenosc}.

\subsection{Nonlinear and nonautonomous stabilization  of Maxwell's equations}

We    consider  Maxwell's equations
in $\Omega \subset  \R^3$ with  a smooth boundary
with either  a nonlinear and nonautonomous internal feedback or a nonlinear and nonautonomous boundary  feedback. To the best of our knowledge, the analysis of Maxwell's system with nonautonomous and nonlinear  damping has not been analyzed.

To clarify the presentation, we treat these two cases separately.

\subsubsection{Interior damping\label{ssMaxint}}
Here we consider the problem
\be\label{max}
\left\{
\begin{tabular}{lllll}
& $\varepsilon \frac{\partial E}{\partial t}
-\rot H+\sigma \sum_{j=1}^J \alpha_j(t,\cdot) g_j(E)=0 $ in $Q$,\vspace{2mm}
\\
& $\mu \frac{\partial H}{\partial t}
+\rot E=0 $ in $Q$,
\\
& $\ddiv(\mu H)=0$ in $Q$,\\
& $E\times \nu=0, H\cdot\nu=0 $
on $\Sigma :=\Gamma \times ]0,+\infty[$,\\
& $E(0) =  E_0,$
$H(0) =  H_0$ in $\Omega $.
\end{tabular}
\right.
\ee
As usual $\eps$ and $\mu$ are real, positive  functions of class  $C^1(\bar\Om)$, while $\sigma$ is a non-negative function that belongs to $L^\infty(\Omega)$ satisfying \eqref{deltacond} 
for some open and non empty subset $\mathcal{O}$ of the the support $X_\sigma$
 of $\sigma$.  For all $j=1,\cdots, J$, the functions $\alpha_j$ and $g_j$ satisfy the assumptions of subsection \ref{ssgeneralhypo} with $U_j=L^2(X_j)^3$, $X_j$ being an open 
 subset of $\Omega$ such that \eqref{partition} holds.

The stability of this system was studied in \cite{Nic-Pignotti:05,nicaise:amo:06,Phung}
with a linear and autonomous feedback $g_j(E)=E$ and $\alpha_j=1$, where some  exponential decay
results were obtained under some constraints on $\eps,\mu$  and $\sigma$. The nonlinear and autonomous case was treated in \cite{nicaise:rendiconti:03}.

Contrary to the above examples this system
is not a second order system
but   it
enters in the setting of (\ref{a6}) once we set
\beqs
\caH&=&L^2(\Omega)^3\times \hat J(\Omega,\mu),\\
\hat J(\Omega,\mu)&=&\{H\in L^2(\Omega)^3: \div(\mu H)=0 \hbox{ in } \Omega,
H\cdot\nu=0
\hbox{ on } \Gamma \},
\\
((E,H),(E',H'))_\caH&=&\int_\Omega (\epsilon E\cdot E'+\mu H\cdot H')\,dx, \forall (E,H),(E',H') \in \caH,\\
\caV&=& V\times \hat J(\Omega,\mu), 
\\
V&=&\{E\in L^2(\Omega)^3: \rot E\in L^2(\Omega)^3, E\times \nu = 0
\hbox{ on } \Gamma \},
\\
\langle A_1(E,H),(E',H')\rangle&=&\int_\Omega (\rot E\cdot H'-H\cdot \rot E')\,dx,
\forall (E,H),(E',H') \in \caV,\\
\langle B(t)(E,H),(E',H')\rangle&=& \sum_{j=1}^J   \int_\Omega \alpha_j(t,x)g_j(E)\cdot E'  \,dx, \forall (E,H),(E',H') \in \caV.
\eeqs

As $\caH$ is continuously embedded into $U=L^2(\Omega)^3\times \{0\}$ (with $I_U (E,H)^\top=(E, 0)^\top$), $\caA_B(t)=A_1+B(t)$ satisfies \eqref{Dt=D}. Furthermore, one readily checks (as in  \cite[\S 3]{ELN}) that $A_1+B(t)$ is maximal monotone
for the inner product $(\cdot,\cdot)_\caH$,   since the bilinear form
$$
\int_\Omega (\mu^{-1}\rot E\cdot \rot E'+\epsilon E\cdot   E')\,dx
$$
is clearly coercive on $V$. Hence, by Theorem  \ref{rkB_2(t)} system (\ref{max})   has a unique solution 
$(E,H)^\top$ of (\ref{max})   satisfying \eqref{regsolutiona6}.

As before  $\pm A_1-I_U$ generates an exponentially stable semigroup in $\caH$ if and only if  system \eqref{max} with a linear and autonomous feedback is  exponentially stable.
As
Theorems 5.1 and 5.5 of \cite{Phung} (resp. Theorem 4.1 of \cite{Nic-Pignotti:05}
and  Remark 5.2 of \cite{nicaise:amo:06})
imply that such an exponential stability holds if $\varepsilon$ and $\mu$ are constant (resp. sufficiently smooth) and under some conditions on $\mathcal{O}$,  
we may conclude  some decays of the solution of \eqref{max} in the setting of one of these papers,
as soon as 
$g_j$ and
$\alpha_j$ satisfy the additional assumptions   from Theorem \ref{tstabobsgen}, \ref{tstabobsgencroissant} or \ref{tstabobsgenosc}.

\subsubsection{Boundary damping\label{ssMaxbdy}}

Let us  go on with     Maxwell's equations   with a nonlinear and non autonous boundary feedback
\be\label{maxbound}
\left\{
\begin{tabular}{lllll}
& $\varepsilon \frac{\partial E}{\partial t}
-\rot H=0 $ in $Q:=\Gamma \times ]0,+\infty[$,\vspace{2mm}
\\
& $\mu \frac{\partial H}{\partial t}
+\rot E=0 $ in $Q$,
\\
& $\ddiv(\varepsilon E)=\ddiv(\mu H)=0$ in $Q$,\\
& $H\times \nu +\alpha(t) g(E\times \nu)\times \nu=0 $
on $\Sigma :=\Gamma \times ]0,+\infty[$,\\
& $E(0) =  E_0,$
$H(0) =  H_0$ in $\Omega $,
\end{tabular}
\right.
\ee
where the functions $\alpha$ and $g$ satisfy the assumptions of subsection \ref{ssgeneralhypo} with $U=L^2(\Gamma)^3$
and $\alpha\in C^1([0,\infty; (0,\infty))$ is such that
  $\alpha'$ is   locally Lipschitz.

The autonomous case
was  studied in \cite{BH,ELN,Kapi,KPan,Komornikbook,Elleretco2002,NicaisePignotti_03,Phung}, where
different decay rates  are avalaible under different conditions on $\epsilon,\mu$ and $\Gamma$
and appropriated assumptions on $g$.

Let us now   show that \eqref{maxbound} enters in the framework of subsection \ref{abunbdcase} if we take (see \cite[\S 2]{Elleretco2002})
\beqs
\caH&=&J(\Omega,\varepsilon)\times  J(\Omega,\mu),\\
J(\Omega,\mu)&=&\{H\in L^2(\Omega)^3: \div(\mu H)=0 \hbox{ in } \Omega\},
\\
((E,H),(E',H'))_\caH&=&\int_\Omega (\epsilon E\cdot E'+\mu H\cdot H')\,dx,\\
\caV&=& V\times J(\Omega,\mu), \forall (E,H),(E',H') \in \caH,
\\
V&=&\{E\in J(\Omega,\varepsilon): \rot E\in L^2(\Omega)^3, E\times \nu \in L^2(\Gamma)^3 \},
\\
U&=& L^2(\Gamma)^3,
\\
\langle A_1(E,H),(E',H')\rangle&=&\int_\Omega (\rot E\cdot H'-H\cdot \rot E')\,dx,
\forall (E,H),(E',H') \in \caV,\\
B(t)&=& \alpha(t)  B_1,\\
\langle B_1(E,H),(E',H')\rangle&=&   \int_\Gamma  g(E\times \nu)\cdot (E'\times \nu)  \,d\sigma(x), \forall (E,H),(E',H') \in \caV.
\eeqs

Note first that $B(t)$ is well-defined with the embedding $I_U(E,H)^\top= E\times \nu$, while by its definition $A_1$ directly satisfies \eqref{h4}. Hence, all assumptions of subsection \ref{ssgeneralhypo} are satisfied. Now in order to apply Theorem \ref{texistencebubdB}, for all $t\geq 0$, we introduce the  bounded linear operators
$D(t)$ and $\tilde D$ from $\caH$ into itself by
\[
D(t)(E,H)^\top=(E, \alpha(t)^{-1} H)^\top,\quad \tilde D(t)(E,H)^\top=(\alpha(t) E,  H)^\top
\]
that, due to the assumption on $\alpha$, satisfy the requested regularity assumptions
and the condition \eqref{equivnorm}  from Theorem \ref{texistencebubdB}.
Furthermore, simple calculations shows that \eqref{commutateur} holds. As Lemma 2.3 of 
 \cite{Elleretco2002} guarantees that the domain of $A_1+B_1$ is dense in $\caH$ and 
 Lemma 2.3 of 
 \cite{Elleretco2002} shows that $A_1+B_1$ is maximal monotone in $\caH$, we can apply Theorem \ref{texistencebubdB} to obtain the well posedness of problem \eqref{maxbound}.
 
 Here again  $\pm A_1-I_U$ generates an exponentially stable semigroup in $\caH$ if and only if  system \eqref{maxbound} with a linear and autonomous feedback is  exponentially stable. Such a stability property was obtained in many papers, let us quote
 \cite{Kapi,KPan,Komornikbook,Phung}.  
 Hence, if  system \eqref{maxbound} with a linear and autonomous feedback is  exponentially stable
 and if additionally
$\alpha$  and $g$   satisfy the additional assumptions   from Theorem \ref{tstabobsgen}, \ref{tstabobsgencroissant} or \ref{tstabobsgenosc}, we may conclude some decays of the solution of \eqref{maxbound}.

        \protect\bibliographystyle{abbrv}
    \protect\bibliographystyle{alpha}
    \bibliography{/Users/Serge/Documents/Serge/Desktop/Biblio/control,/Users/Serge/Documents/Serge/Desktop/Biblio/valein,/Users/Serge/Documents/Serge/Desktop/Biblio/kunert,/Users/Serge/Documents/Serge/Desktop/Biblio/est,/Users/Serge/Documents/Serge/Desktop/Biblio/femaj,/Users/Serge/Documents/Serge/Desktop/Biblio/mgnet,/Users/Serge/Documents/Serge/Desktop/Biblio/dg,/Users/Serge/Documents/Serge/Desktop/Biblio/bib,/Users/Serge/Documents/Serge/Desktop/Biblio/maxwell,/Users/Serge/Documents/Serge/Desktop/Biblio/bibmix,/Users/Serge/Documents/Serge/Desktop/Biblio/cochez,/Users/Serge/Documents/Serge/Desktop/Biblio/soualem,/Users/Serge/Documents/Serge/Desktop/Biblio/nic}

\newcommand{\noopsort}[1]{}\def\ocirc#1{\ifmmode\setbox0=\hbox{$#1$}\dimen0=\ht0
  \advance\dimen0 by1pt\rlap{\hbox to\wd0{\hss\raise\dimen0
  \hbox{\hskip.2em$\scriptscriptstyle\circ$}\hss}}#1\else {\accent"17 #1}\fi}
  \def\cprime{$'$} \def\polhk#1{\setbox0=\hbox{#1}{\ooalign{\hidewidth
  \lower1.5ex\hbox{`}\hidewidth\crcr\unhbox0}}} \def\cprime{$'$}
\begin{thebibliography}{10}

\bibitem{AlabauKomornik}
F.~Alabau and V.~Komornik.
\newblock Boundary observability, controllability, and stabilization of linear
  elastodynamic systems.
\newblock {\em SIAM J. Control Optim.}, 37(2):521--542 (electronic), 1999.

\bibitem{ammari:01b}
K.~Ammari, A.~Henrot, and M.~Tucsnak.
\newblock Asymptotic behaviour of the solutions and optimal location of the
  actuator for the pointwise stabilization of a string.
\newblock {\em Asymptot. Anal.}, 28(3-4):215--240, 2001.

\bibitem{ammari:01}
K.~Ammari and M.~Tucsnak.
\newblock Stabilization of second order evolution equations by a class of
  unbounded feedbacks.
\newblock {\em ESAIM Control Optim. Calc. Var.}, 6:361--386, 2001.

\bibitem{bardos}
C.~Bardos, G.~Lebeau, and J.~Rauch.
\newblock Sharp sufficient conditions for the observation, control, and
  stabilization of waves from the boundary.
\newblock {\em SIAM J. Control Optim.}, 30(5):1024--1065, 1992.

\bibitem{BH}
H.~Barucq and B.~Hanouzet.
\newblock Etude asymptotique du syst\`eme de {M}axwell avec la condition aux
  limites absorbante de {S}ilver-{M}\"uller {II}.
\newblock {\em C. R. Acad. Sci. Paris S\'er. I}, 316:1019--1024, 1993.

\bibitem{BchatniaDaoulatli:13}
A.~Bchatnia and M.~Daoulatli.
\newblock Behavior of the energy for {L}am\'{e} systems in bounded domains with
  nonlinear damping and external force.
\newblock {\em Electron. J. Differential Equations}, pages No. 01, 17, 2013.

\bibitem{Bellassoued:08}
M.~Bellassoued.
\newblock Energy decay for the elastic wave equation with a local
  time-dependent nonlinear damping.
\newblock {\em Acta Math. Sin. (Engl. Ser.)}, 24(7):1175--1192, 2008.

\bibitem{Bey:03}
R.~Bey, A.~Heminna, and J.-P. Loh{\'e}ac.
\newblock Boundary stabilization of the linear elastodynamic system by a
  {L}yapunov-type method.
\newblock {\em Rev. Mat. Complut.}, 16(2):417--441, 2003.

\bibitem{Browder:68}
F.~E. Browder.
\newblock Nonlinear maximal monotone operators in {B}anach space.
\newblock {\em Math. Ann.}, 175:89--113, 1968.

\bibitem{Chen_G:79b}
G.~Chen.
\newblock Control and stabilization for the wave equation in a bounded domain.
\newblock {\em SIAM J. Control Optim.}, 17(1):66--81, 1979.

\bibitem{Chen_G:79a}
G.~Chen.
\newblock Energy decay estimates and exact boundary value controllability for
  the wave equation in a bounded domain.
\newblock {\em J. Math. Pures Appl. (9)}, 58(3):249--273, 1979.

\bibitem{Crandal72}
M.~G. Crandall and A.~Pazy.
\newblock Nonlinear evolution equations in {B}anach spaces.
\newblock {\em Israel J. Math.}, 11:57--94, 1972.

\bibitem{daLuz-PerlaMenzala:11}
C.~R. da~Luz and G.~A. Perla~Menzala.
\newblock Uniform decay rates of coupled anisotropic elastodynamic/{M}axwell
  equations with nonlinear damping.
\newblock {\em Port. Math.}, 68(2):205--238, 2011.

\bibitem{Daoutli:11}
M.~Daoulatli.
\newblock Rates of decay for the wave systems with time dependent damping.
\newblock {\em Discrete Contin. Dyn. Syst.}, 31(2):407--443, 2011.

\bibitem{ELN}
M.~Eller, J.~Lagnese, and S.~Nicaise.
\newblock Decay rates for solutions of a {M}axwell system with nonlinear
  boundary damping.
\newblock {\em Comp. and Appl. Math.}, 21:135--165, 2002.

\bibitem{Evans77}
L.~C. Evans.
\newblock Nonlinear evolution equations in an arbitrary banach space.
\newblock {\em Israel J. Math.}, 26:1--42, 1977.

\bibitem{Guesmia98a}
A.~Guesmia.
\newblock Existence globale et stabilisation interne non lin\'eaire d'un
  syst\`eme de {P}etrovsky.
\newblock {\em Bull. Belg. Math. Soc. Simon Stevin}, 5(4):583--594, 1998.

\bibitem{guesmia98b}
A.~Guesmia.
\newblock Observability, controllability and boundary stabilization of some
  linear elasticity systems.
\newblock {\em Acta Sci. Math. (Szeged)}, 64(1-2):109--119, 1998.

\bibitem{Guesmia99b}
A.~Guesmia.
\newblock Existence globale et stabilisation fronti\`ere non lin\'eaire d'un
  syst\`eme d'\'elasticit\'e.
\newblock {\em Portugal. Math.}, 56(3):361--379, 1999.

\bibitem{guesmia:00}
A.~Guesmia.
\newblock On the decay estimates for elasticity systems with some localized
  dissipations.
\newblock {\em Asymptot. Anal.}, 22(1):1--13, 2000.

\bibitem{haraux:89a}
A.~Haraux.
\newblock Une remarque sur la stabilisation de certains syst\`emes du
  deuxi\`eme ordre en temps.
\newblock {\em Portugal. Math.}, 46(3):245--258, 1989.

\bibitem{Horn}
M.~A. Horn.
\newblock Implications of sharp trace regularity results on boundary
  stabilization of the system of linear elasticity.
\newblock {\em J. Math. Anal. Appl.}, 223(1):126--150, 1998.

\bibitem{Jiao-Xiao:15}
Z.~Jiao and T.-J. Xiao.
\newblock Convergence and speed estimates for semilinear wave systems with
  nonautonomous damping.
\newblock {\em Math. Methods Appl. Sci.}, 39(18):5465--5474, 2016.

\bibitem{Kapi}
B.~V. Kapitanov.
\newblock Stabilization and exact boundary controllability for {M}axwell's
  equations.
\newblock {\em SIAM J. Control Optim.}, 32:408--420, 1994.

\bibitem{KapiRaupp}
B.~V. Kapitanov and M.~A. Raupp.
\newblock Exact boundary controllability in problems of transmission for the
  system of electromagneto-elasticity.
\newblock {\em Math. Meth. Appl. Sci.}, 24:193--207, 2001.

\bibitem{Kato67}
T.~Kato.
\newblock Nonlinear semigroups and evolution equations.
\newblock {\em J. Math. Soc. Japan}, 19:508--520, 1967.

\bibitem{KatoPisa}
T.~Kato.
\newblock {\em Abstract differential equations and nonlinear mixed problems}.
\newblock Lezioni Fermiane. [Fermi Lectures]. Scuola Normale Superiore, Pisa;
  Accademia Nazionale dei Lincei, Rome, 1985.

\bibitem{KatoCIME}
T.~Kato.
\newblock Linear and quasi-linear equations of evolution of hyperbolic type.
\newblock In {\em Hyperbolicity}, volume~72 of {\em C.I.M.E. Summer Sch.},
  pages 125--191. Springer, Heidelberg, 2011.

\bibitem{Komornik91}
V.~Komornik.
\newblock Rapid boundary stabilization of the wave equation.
\newblock {\em SIAM J. Control Optim.}, 29:197--208, 1991.

\bibitem{Komornik93b}
V.~Komornik.
\newblock On the nonlinear boundary stabilization of the wave equation.
\newblock {\em Chinese Ann. Math. Ser. B}, 14:153--164, 1993.

\bibitem{KPan}
V.~Komornik.
\newblock Boundary stabilization, observation and control of {M}axwell's
  equations.
\newblock {\em PanAm. Math. J.}, 4:47--61, 1994.

\bibitem{Komornik93a}
V.~Komornik.
\newblock {\em Decay estimates for the wave equation with internal damping},
  volume 118 of {\em Int. Ser. Numer. Anal.}, pages 253--266.
\newblock Vorau, 1994.

\bibitem{Komornikbook}
V.~Komornik.
\newblock {\em Exact controllability and stabilization}.
\newblock RAM: Research in Applied Mathematics. Masson, Paris, 1994.
\newblock The multiplier method.

\bibitem{komornik:94}
V.~Komornik.
\newblock On the nonlinear boundary stabilization of {K}irchhoff plates.
\newblock {\em NoDEA Nonlinear Differential Equations Appl.}, 1(4):323--337,
  1994.

\bibitem{KomornikZuazua}
V.~Komornik and E.~Zuazua.
\newblock A direct method for the boundary stabilization of the wave equation.
\newblock {\em J. Math. Pures Appl. (9)}, 69(1):33--54, 1990.

\bibitem{LaSiam83}
J.~Lagnese.
\newblock Boundary stabilization of linear elastodynamic systems.
\newblock {\em SIAM J. Control Optim.}, 21(6):968--984, 1983.

\bibitem{LaJDE83}
J.~Lagnese.
\newblock Decay of solutions of wave equations in a bounded region with
  boundary dissipation.
\newblock {\em J. Differential Equations}, 50(2):163--182, 1983.

\bibitem{Lagnese:88}
J.~E. Lagnese.
\newblock Note on boundary stabilization of wave equations.
\newblock {\em SIAM J. Control Optim.}, 26(5):1250--1256, 1988.

\bibitem{Lasiecka-Tataru:93}
I.~Lasiecka and D.~Tataru.
\newblock Uniform boundary stabilization of semilinear wave equations with
  nonlinear boundary damping.
\newblock {\em Differential Integral Equations}, 6(3):507--533, 1993.

\bibitem{Lasiecka-Triggiani92}
I.~Lasiecka and R.~Triggiani.
\newblock Uniform stabilization of the wave equation with {D}irichlet or
  {N}eumann feedback control without geometrical conditions.
\newblock {\em Appl. Math. Optim.}, 25(2):189--224, 1992.

\bibitem{lebeau:96}
G.~Lebeau.
\newblock \'{E}quation des ondes amorties.
\newblock In {\em Algebraic and geometric methods in mathematical physics
  ({K}aciveli, 1993)}, volume~19 of {\em Math. Phys. Stud.}, pages 73--109.
  Kluwer Acad. Publ., Dordrecht, 1996.

\bibitem{Lin}
C.-Y. Lin.
\newblock Time-dependent nonlinear evolution equations.
\newblock {\em Differential Integral Equations}, 15:257--270, 2002.

\bibitem{Liu}
K.~Liu and Z.~Liu.
\newblock Boundary stabilization of a nonhomogeneous beam with rotatory inertia
  at the tip.
\newblock {\em J. Comput. Appl. Math.}, 114(1):1--10, 2000.
\newblock Control of partial differential equations (Jacksonville, FL, 1998).

\bibitem{Liuzuazua}
W.-J. Liu and E.~Zuazua.
\newblock Decay rates for dissipative wave equations.
\newblock volume~48, pages 61--75. 1999.
\newblock Papers in memory of Ennio De Giorgi (Italian).

\bibitem{Luo-Xiao:20a}
J.-R. Luo and T.-J. Xiao.
\newblock Decay rates for second order evolution equations in {H}ilbert spaces
  with nonlinear time-dependent damping.
\newblock {\em Evol. Equ. Control Theory}, 9(2):359--373, 2020.

\bibitem{Luo-Xiao:20b}
J.-R. Luo and T.-J. Xiao.
\newblock Decay rates for semilinear wave equations with vanishing damping and
  {N}eumann boundary conditions.
\newblock {\em Math. Methods Appl. Sci.}, 44(1):303--314, 2021.

\bibitem{Luo-Xiao:21}
J.-R. Luo and T.-J. Xiao.
\newblock Optimal energy decay rates for abstract second order evolution
  equations with non-autonomous damping.
\newblock {\em ESAIM Control Optim. Calc. Var.}, 27:Paper No. 59, 24, 2021.

\bibitem{martinez:99}
P.~Martinez.
\newblock A new method to obtain decay rate estimates for dissipative systems
  with localized damping.
\newblock {\em Rev. Mat. Comp. Madrid}, 12:251--283, 1999.

\bibitem{martinez:00}
P.~Martinez.
\newblock Precise decay rate estimates for time-dependent dissipative systems.
\newblock {\em Israel J. Math.}, 119:291--324, 2000.

\bibitem{mustapha:15}
M.~I. Mustafa.
\newblock Uniform decay for wave equations with weakly dissipative boundary
  feedback.
\newblock {\em Dyn. Syst.}, 30(2):241--250, 2015.

\bibitem{Nakao:97}
M.~Nakao.
\newblock On the decay of solutions of the wave equation with a local
  time-dependent nonlinear dissipation.
\newblock {\em Adv. Math. Sci. Appl.}, 7(1):317--331, 1997.

\bibitem{nicaise:rendiconti:03}
S.~Nicaise.
\newblock Stability and controllability of an abstract evolution equation of
  hyperbolic type and concrete applications.
\newblock {\em Rendiconti di Matematica {S}erie VII}, 23:83--116, 2003.

\bibitem{nicEES}
S.~Nicaise.
\newblock Stability and controllability of the electromagneto-elastic system.
\newblock {\em Port. Math. (N.S.)}, 60(1):37--70, 2003.

\bibitem{Elleretco2002}
S.~Nicaise, M.~Eller, and J.~E. Lagnese.
\newblock Stabilization of heterogeneous {M}axwell's equations by linear or
  nonlinear boundary feedback.
\newblock {\em Electron. J. Differential Equations}, pages No. 21, 26 pp.
  (electronic), 2002.

\bibitem{NicaisePignotti_03}
S.~Nicaise and C.~Pignotti.
\newblock Boundary stabilization of {M}axwell's equations with space-time
  variable coefficients.
\newblock {\em ESAIM Control Optim. Calc. Var.}, 9:563--578 (electronic), 2003.

\bibitem{Nic-Pignotti:05}
S.~Nicaise and C.~Pignotti.
\newblock Internal stabilization of {M}axwell's equations in heterogeneous
  media.
\newblock {\em Abstr. Appl. Anal.}, (7):791--811, 2005.

\bibitem{nicaise:amo:06}
S.~Nicaise and C.~Pignotti.
\newblock Internal and boundary observability estimates for heterogeneous
  maxwell's system.
\newblock {\em Applied Math. Optim.}, 54:47--70, 2006.

\bibitem{Pazy}
A.~Pazy.
\newblock {\em Semigroups of linear operators and applications to partial
  differential equations}, volume~44 of {\em Applied Math. Sciences}.
\newblock Springer-Verlag, New York, 1983.

\bibitem{Phung}
K.~D. Phung.
\newblock Contr\^ole et stabilisation d'ondes \'electromagn\'etiques.
\newblock {\em ESAIM: COCV}, 5:87--137, 2000.

\bibitem{Pucci-Serrin:96}
P.~Pucci and J.~Serrin.
\newblock Asymptotic stability for nonautonomous dissipative wave systems.
\newblock {\em Comm. Pure Appl. Math.}, 49(2):177--216, 1996.

\bibitem{Tebou:96}
L.~R. Tcheugou\'{e}~Tebou.
\newblock On the stabilization of the wave and linear elasticity equations in
  {$2$}-{D}.
\newblock {\em PanAmer. Math. J.}, 6(1):41--55, 1996.

\bibitem{Triggiani:89}
R.~Triggiani.
\newblock Wave equation on a bounded domain with boundary dissipation: an
  operator approach.
\newblock {\em J. Math. Anal. Appl.}, 137(2):438--461, 1989.

\bibitem{zuazua:cpde:90}
E.~Zuazua.
\newblock Exponential decay for the semilinear wave equation with locally
  distributed damping.
\newblock {\em Comm. in PDE}, 15:205--235, 1990.

\bibitem{zuazua:sicon:90}
E.~Zuazua.
\newblock Uniform stabilization of the wave equation by nonlinear boundary
  feedback.
\newblock {\em SIAM J. Control Optim.}, 28:466--477, 1990.

\end{thebibliography}

\end{document}